\documentclass[11pt,reqno]{amsart}
\newtheorem{theorem}{Theorem}[section]
\newtheorem{proposition}[theorem]{Proposition}
\newtheorem{corollary}[theorem]{Corollary}
\newtheorem{remark}[theorem]{Remark}
\newtheorem{lemma}[theorem]{Lemma}
\newtheorem{example}[theorem]{Example}
\newtheorem{definition}[theorem]{Definition}

\numberwithin{equation}{section}
\usepackage{amssymb}
\usepackage{graphicx}
\usepackage{amsmath}
\usepackage{enumerate}
\usepackage{color}
\usepackage[all]{xy}
\usepackage[latin1]{inputenc} 
\usepackage{graphicx} 
\usepackage{mathtools}
\begin{document}

\title[Lipeng Luo \textsuperscript{1} and Sania Asif\textsuperscript{2}]{On some derivations of Lie conformal superalgebras}
\author{Lipeng Luo \textsuperscript{1} and Sania Asif\textsuperscript{2}}

\address{\textsuperscript{1}Department of Mathematics, Tongji University, Shanghai, 200092, PR China.}
\address{\textsuperscript{2}School of Mathematics and Statistics, Nanjing University of Information Science and Technology, Nanjing, Jiangsu Province, 210044, PR China.}

\email{\textsuperscript{1}luolipeng@tongji.edu.cn}
\email{\textsuperscript{2}200036@nuist.edu.cn}




\date{\today}

\begin{abstract}
Let $\mathcal{R}$ be a Lie conformal superalgebra. In this paper, we first investigate the conformal derivation algebra $CDer(\mathcal{R})$, the conformal triple derivation algebra $CTDer(\mathcal{R})$, and the generalized conformal triple derivation algebra $GCTDer(\mathcal{R})$. Moreover, we determine the connection of these derivation algebras. Next, we give a complete classification of the (generalized) conformal triple derivation algebra on all finite simple Lie conformal superalgebras. More specifically, $CTDer(\mathcal{R})=CDer(\mathcal{R})$, where $\mathcal{R}$ is a finite simple Lie conformal superalgebra, but for $GCTDer(\mathcal{R})$, we obtain a conclusion that is closely related to $CDer(\mathcal{R})$. 
Furthermore, we evaluate the $(\varPhi, \varPsi)$-Lie triple derivations on Lie conformal superalgebra, where $\varPhi$ and $\varPsi$ are associated automorphism of $\phi_{x}\in gc(\mathcal R)$. We evaluated some fundamental properties of $(\varPhi, \varPsi)$- Lie triple derivations.
 Later, we introduce the definition of $(A, B, C, D)$-derivation on Lie conformal superalgebra. We obtain the relationships between the generalized conformal triple derivations and the conformal $(A, 
B, C, D)$-derivations on Lie conformal superalgebra.
Finally, we have presented the triple homomorphism of Lie conformal superalgebras.
\end{abstract}

\footnote{The second author is the corresponding author: 200036@nuist.edu.cn(Asif S.)}
\maketitle

\section{Introduction}\label{intro}
The structure of Lie conformal superalgebra was first studied by Fattori in \cite{Fattori-Kac}, where representations and the derivations of Lie conformal superalgebra were studied. Recently, the theory of generalized derivations of  Lie conformal algebras are considered in \cite{Fan-Hong-Su} and the  generalized derivations of Lie algebras and extension to the theory are given in \cite{Leger-Lucks}. Lie conformal superalgebra has a special role in the operator product of the chiral fields (2D quantum field theory). Moreover, there is a close connection between Lie conformal superalgebras and the concept of a formal distribution Lie superalgebra $(\mathcal{G}, F)$, i.e., a Lie superalgebra $\mathcal{G}$ is spanned by the coefficients of a family $F$ of mutually local formal distributions, see \cite{Boyallian-Kac-Liberati-Rudakov, Chang, Fattori-Kac} for details.

In the past few years, a remarkable study has been done on the  derivations and generalized derivations of Lie algebra, Lie superalgebra, Lie color algebra, Lie conformal algebra, Lie triple system, $n$-Hom Nambu Lie algebra, and Hom-Lie triple system  in \cite{Chen-Ma-Ni, Hopkins, Leger-Lucks, Zhang-Zhang, Zhou-Chen, Zhou-Chen-Ma-2, Zhou-Chen-Ma-3, Zhao-Chen-Yuan}. It is prominent that this research on derivations and generalized derivations algebra of Lie algebra plays a significant role in the development of the structural theory of Lie algebras. Leger and Lucks investigated the generalized derivations of Lie algebras and their subalgebras in \cite{Leger-Lucks}. There are some generalizations of the derivation of Lie algebra, such as the biderivations of the Lie algebra, Schr\"{o}dinger-Virasoro Lie algebra and Block Lie algebras, that were investigated by various authors in \cite{Bresar-Zhao, Liu-Guo-Zhao, Wang-Yu}. Lie triple derivation of Lie algebra as a generalization of derivation was first introduced by M\"{u}ller in \cite{Muller}. Triple derivations of Lie algebra resemble the triple derivations of Jordan algebra and associative algebra.  It is not hard to see that, Lie derivations are in fact  Lie triple derivations, but converse to this statement may or may not be true. After M\"{u}ller, many authors have studied Lie triple derivations from various points of view and determined many interesting results, i.e., study on Lie triple derivations of nest algebra, quaternion algebra, triangular algebra, Lie color algebra, TUHF algebra, Lie algebra of strictly upper triangular matrix over a commutative ring, $gl(n,\mathcal{R})$ and perfect Lie (super)algebras, see  \cite{Asif-Wu, Ji-Wang, Lu, Li-Wang, Li-Wang-Guo, Wu-Asif-Munir, Xiao-Wei, Zhou,  Zhou-Chen-Ma-1} for details. There are still some other triple derivations on associative algebras, such as Jordan triple derivations and associative triple derivations in \cite{Beidar-Chebotar, Herstein}.

Lately, the $(\alpha, \beta, \gamma)$-derivations of Lie algebras, Lie superalgebras, and Lie conformal superalgebra are discussed in \cite{Feng-Zhao-Chen, Novotny-Hrivnak, Zheng-Zhang}. In this paper, we focused on studying triple derivations, generalized triple derivations, centroid derivation, and central triple derivations of Lie conformal superalgebra. Moreover, we also study the triple derivations of finite Lie conformal superalgebra. We evaluate that, for the finite Lie conformal superalgebra, corresponding triple derivations algebra coincides with derivations algebra. Furthermore, by following \cite{Feng-Zhao-Chen}, we evaluate $(\varPhi,\varPsi)$ triple derivations of Lie conformal superalgebra and give many interesting properties of these derivations. We also aim to study $(A,B,C,D)$-derivations of Lie conformal superalgebra. We show that, by applying different constraints, $(A, B, C, D)$-derivations show numerous properties. 

This paper is arranged as follows. In Section $2$, we present a few fundamental notations, definitions,  and related results about Lie conformal superalgebras. In Section $3$, firstly we focus on establishing the notation of (generalized) conformal triple derivations, centroid derivations, and central triple derivations on Lie conformal superalgebras. Then we investigate the Lie conformal superalgebra structure and some properties of these triple derivations. In Section $4$, we investigate the (generalized) conformal triple derivations of all finite Lie conformal superalgebras and give the complete classification for these derivations. In Section $5$, we first introduce the notion of $(\varPhi, \varPsi)$-derivations of Lie conformal superalgebra and then give the fundamental properties of this derivation system with the help of interesting propositions. Furthermore, we present the definition and classification of the $(A, B, C, D)$-derivations of Lie conformal superalgebra and find out many important identities. In Section $6$, we first introduce the definition of triple homomorphism of a Lie conformal superalgebra. Then, we give the classification of all finite simple Lie conformal superalgebras. 
\par Throughout this paper, $\mathbb{C}$ represents the set of complex numbers. On the other hand, all  tensor products and vector spaces are considered over $\mathbb{C}$. Without any  uncertainty, we abbreviate $\otimes_{\mathbb{C}}$ by $\otimes$.

\section{preliminaries}

 In this section, we recall some basic definitions, notations and related results about Lie conformal superalgebras for later use. For a detailed description, one can refer to \cite{D'Andrea-Kac, KacV, Zhao-Chen-Yuan}.
 
Consider a $Z_{2}$-graded linear space with a direct sum $V= V_{\overline{0}}\oplus V_{\overline{1}}$ denoted by $V$.  $V$ is a  superspace, having the homogeneous elements with the parity $j$, denoted by $V_j$, where $ j\in\{\overline{0}, \overline{1}\}$. For example, if $j$ is the parity of $a$, then we write  $|a|=j$. In this paper $|a|$ denotes the parity of homogeneous element $a$ and this expression can extend to other elements by linearity. 
\begin{definition}\label{def2.1}
\begin{em}
A \emph {Lie conformal superalgebra} $\mathcal{R}$ is a $Z_{2}$-graded $\mathbb{C}[\partial]$-module endowed with a $\mathbb{C}$-linear map from $\mathcal{R}\otimes\mathcal{R} \to \mathbb{C}[\lambda]\otimes\mathcal{R}$, $a\otimes b \mapsto [a_\lambda b]$, called the $\lambda$-bracket, satisfying the following axioms:
\begin{align}
    [\partial a_\lambda b]&=-\lambda[a_\lambda b],\quad  [a_\lambda \partial b]=(\partial+\lambda)[a_\lambda b] \quad (conformal \  sesquilinearity),\\
    {}[a_\lambda b]&=-(-1)^{|a||b|}[b_{-\lambda-\partial} a]\quad (skew\text{-}symmetry),\\
    {}[a_\lambda [b_\mu c]]&=[[a_\lambda b]_{\lambda+\mu} c]+(-1)^{|a||b|}[ b_\mu [a_\lambda c]]\quad (Jacobi\  identity),
\end{align}
for $a,b,c \in\mathcal{R}$.
\end{em}		
\end{definition}

 As shown in \cite{KacV}, defining the $\lambda$-bracket as follows:
 \begin{align}
 [a_{\lambda}b]=\sum_{n\in \mathbb{Z^+}}\frac{\lambda^n}{n!}a_{(n)}b,\quad \forall a,b \in \mathcal{R}.
 \end{align}
 The equivalent definition of Lie conformal superalgebra can be written as follows.
\begin{definition}\label{def2.2}
\begin{em}
A \emph {Lie conformal superalgebra} $\mathcal{R}$ is a $Z_{2}$-graded $\mathbb{C}[\partial]$-module endowed with a $\mathbb{C}$-bilinear map from $\mathcal{R}\otimes\mathcal{R} \to \mathcal{R}$, called $n$-product, for any $n\in \mathbb{Z^+}$, satisfying ${\mathcal{R}_{\alpha}}_{(n)}{\mathcal{R}_\beta}\subseteq {\mathcal{R}_{\alpha+\beta}}$, where $\alpha, \beta \in \{\overline{0}, \overline{1}\}$ and the following axioms: for any $a,b\in\mathcal{R}$, and $m,n\in \mathbb{Z^+}$, 
\begin{enumerate}[(C1)]
\item $a_{(n)}b=0$ for $n\gg 0$,
\item $(\partial a)_{(n)} b=-na_{(n-1)} b$, and $a_{(n)}(\partial b)=\partial (a_{(n)}b)+na_{(n-1)} b$,
\item $a_{(n)}b=-(-1)^{|a||b|}\sum_{j\in\mathbb{Z^+}}(-1)^{n+j}\frac{\partial^j(b_{(n+j)}a)}{j!}$,
\item $a_{(m)}(b_{(n)}c)=\sum_{j=0}^{m}\binom m j (a_{(j)}b)_{(m+n-j)}c+(-1)^{|a||b|}b_{(n)}(a_{(m)}c)$.
 \end{enumerate}
\end{em}		
\end{definition}

A \emph{subalgebra} $\mathcal{S}$ of $\mathcal{R}$ is a $\mathbb{C}[\partial]$-submodule of $\mathcal{R}$ such that $\mathcal{S}_{(n)}\mathcal{S}\subseteq\mathcal{S}$ for any $n\in \mathbb{Z^+}$. An \emph{ideal} $\mathcal{I}$ of $\mathcal{R}$ is a $\mathbb{C}[\partial]$-submodule of $\mathcal{R}$ such that $\mathcal{R}_{(n)}\mathcal{I}\subseteq\mathcal{I}$ for any $n\in \mathbb{Z^+}$. Moreover, due to the skew\text{-}symmetry, any left or right ideal is actually a two-side ideal. A Lie conformal superalgebra $\mathcal{R}$ is called \emph {finite} if $\mathcal{R}$ is finitely generated as a $\mathbb{C}[\partial]$-module. The \emph {rank} of a Lie conformal superalgebra $\mathcal{R}$, denoted by rank($\mathcal{R}$), is its rank as a $\mathbb{C}[\partial]$-module. A Lie conformal superalgebra $\mathcal{R}$ is \emph{simple} if it has no non-trivial ideals and it is not abelian. Some examples of Lie conformal superalgebras are given as follows.

\begin{example}\label{ex2.2}
\begin{em}
Let $\mathcal{R}=\mathbb{C}[\partial]L\oplus\mathbb{C}[\partial]W$ be a free $Z_{2}$-graded $\mathbb{C}[\partial]$-module. Define 
\begin{align*}
[L_\lambda L]=(\partial+2\lambda)L, \quad [L_\lambda W]=(\partial+\frac{3}{2}\lambda)W, \quad [W_\lambda W]=0,
\end{align*}
where $\mathcal{R}_{\overline{0}}=\mathbb{C}[\partial]L$ and $\mathcal{R}_{\overline{1}}=\mathbb{C}[\partial]W$. Then $\mathcal{R}$ is a Lie conformal superalgebra.
\end{em}	
\end{example}

\begin{example}\label{ex2.3}
\begin{em}
(Neveu-Schwarz Lie conformal superalgebra)Let $\mathcal{NS}=\mathbb{C}[\partial]L\oplus\mathbb{C}[\partial]G$ be a free $Z_{2}$-graded $\mathbb{C}[\partial]$-module. Define 
\begin{align*}
[L_\lambda L]=(\partial+2\lambda)L, \quad [L_\lambda G]=(\partial+\frac{3}{2}\lambda)G, \quad [G_\lambda G]=2L,
\end{align*}
where $\mathcal{NS}_{\overline{0}}=\mathbb{C}[\partial]L$ and $\mathcal{NS}_{\overline{1}}=\mathbb{C}[\partial]G$. Then $\mathcal{NS}$ is a Lie conformal superalgebra. We call it Neveu-Schwarz Lie conformal superalgebra.
\end{em}	
\end{example}

\begin{definition}\label{def2.3}
\begin{em}
A homogeneous $\mathbb{C}[\partial]$-module homomorphism $f$: $\mathcal{A} \to \mathcal{B}$ is called a \emph {homomorphism} between Lie conformal superalgebras $\mathcal{A}$ and $\mathcal{B}$ if it satisfies $f([a_\lambda b])=[f(a)_\lambda f(b)]$ for any $a,b\in\mathcal{A}$.
\end{em}		
\end{definition}

\begin{example}\label{ex2.4}
\begin{em}
Let $\mathcal{G}=\mathcal{G}_{\overline{0}}\oplus\mathcal{G}_{\overline{1}}$ be a complex Lie superalgebra. Set $(Cur\mathcal{G})_\alpha:=\mathbb{C}[\partial]\otimes\mathcal{G}_\alpha$ as a free $\mathbb{C}[\partial]$-module. Then $Cur\mathcal{G}=(Cur\mathcal{G})_{\overline{0}}\oplus(Cur\mathcal{G})_{\overline{1}}$ is a Lie conformal superalgebra, called \emph{current Lie conformal superalgebra}, with $\lambda$-bracket given by 
\begin{align}
[(f(\partial)\otimes a)_\lambda(g(\partial)\otimes b)]:=f(-\lambda)g(\partial+\lambda)\otimes[a,b], \quad \forall a,b \in \mathcal{R}.
\end{align}
\end{em}	
\end{example}
It was shown in \cite{Fattori-Kac} that all finite simple Lie conformal superalgebras were determined. In this paper, we will classify conformal triple derivations and conformal triple homomorphisms of all finite simple Lie conformal superalgebras.

\begin{definition}\label{def2.4}
\begin{em}
A \emph {conformal linear map of degree $\alpha$} between $Z_{2}$-graded $\mathbb{C}[\partial]$-modules $\mathcal{A}$ and $\mathcal{B}$ is a $\mathbb{C}$-linear map $\phi_{\lambda}$: $\mathcal{A} \to \mathbb{C}[\lambda]\otimes\mathcal{B}$, satisfying the following axioms:
\begin{align*}
\phi_{\lambda}(\partial a)=(\partial+\lambda)\phi_{\lambda}(a), \quad \forall  a\in\mathcal{A}, \quad and\quad \phi_{\lambda}(\mathcal{A}_\beta)\subseteq\mathcal{B}_{\alpha+\beta}[\lambda], \quad \forall \alpha, \beta \in \{\overline{0}, \overline{1}\}.
 \end{align*}
\end{em}		
\end{definition}

Obviously, $\phi$ is conformal linear, which does not depend on the choice of the indeterminate variable $\lambda$. Throughout what follows, we will also write $\phi_{\lambda}$ instead of $\phi$ to emphasize the dependence of $\phi$ on $\lambda$. The vector space of all conformal linear maps of degree $\alpha$ from $\mathcal{A}$ to $\mathcal{B}$ is denoted by $Chom(\mathcal{A},\mathcal{B})_\alpha$. Then $Chom(\mathcal{A},\mathcal{B})=Chom(\mathcal{A},\mathcal{B})_{\overline{0}}\oplus Chom(\mathcal{A},\mathcal{B})_{\overline{1}}$, which can be made into a $Z_{2}$-graded $\mathbb{C}[\partial]$-module via $(\partial \phi)_{\lambda}(a)=-\lambda\phi_{\lambda}(a), \forall  a\in\mathcal{A}$. For convenience, we write $Cend{\mathcal{A}}$ for $Chom(\mathcal{A},\mathcal{A})$. 

\begin{definition}\label{def2.5}
\begin{em}
Let $\mathcal{R}$ be a Lie conformal superalgebra. A conformal linear map $\phi_{\lambda}\in (Cend{\mathcal{R}})_\alpha$ is a \emph {conformal derivation of degree $\alpha$} of $\mathcal{R}$ if 
\begin{align*}
\phi_{\lambda}([a_{\mu}b])=[\phi_{\lambda}(a)_{\lambda+\mu}b]+(-1)^{|\phi||a|}[a_{\mu}\phi_{\lambda}(b)],\quad \forall a,b \in \mathcal{R}.
\end{align*}
\end{em}		
\end{definition}

By the Jacobi identity, it is not difficult to check that for every $a \in \mathcal{R}$, the map $ad a _{\lambda}$ which is defined by $ad a _{\lambda} b=[a_{\lambda} b]$ for any $b \in \mathcal{R}$, is a conformal derivation of $\mathcal{R}$. All conformal derivations of this kind are called \emph{inner conformal derivations}, Denote by $CDer(\mathcal{R})$ and $CInn(\mathcal{R})$ the vector spaces of all conformal derivations and inner conformal derivations of $\mathcal{R}$ respectively. Moreover, we have $CDer(\mathcal{R})=CDer(\mathcal{R})_{\overline{0}}\oplus CDer(\mathcal{R})_{\overline{1}}$ and $CInn(\mathcal{R})=CInn(\mathcal{R})_{\overline{0}}\oplus CInn(\mathcal{R})_{\overline{1}}$

As shown in \cite{Zhao-Chen-Yuan}, we have an important result listed below.
\begin{proposition}\label{ex2.6}
Every conformal derivation of the Neveu-Schwarz Lie conformal superalgebra is an inner conformal derivation.	
\end{proposition}

\begin{definition}\label{def2.7}
\begin{em}
Let $\mathcal{R}$ be a Lie conformal superalgebra. The $\lambda$-bracket on $Cend(\mathcal{R})$ given by 
\begin{align*}
[\phi_{\lambda}\psi]_{\mu}a=\phi_{\lambda}(\psi_{\mu-\lambda}a)-(-1)^{|\phi||\psi|}\psi_{\mu-\lambda}(\phi_{\lambda}a),\quad \forall a \in \mathcal{R},
\end{align*}
defines a Lie conformal superalgebra structure on $Cend(\mathcal{R})$. This is called the \emph{general Lie conformal superalgebra} on $\mathcal{R}$ which is denoted by $gc(\mathcal{R})$. It is obvious that $gc(\mathcal{R})=gc(\mathcal{R})_{\overline{0}}\oplus gc(\mathcal{R})_{\overline{1}}$.
\end{em}		
\end{definition}
Inspired by the above definition, it is easy to check that $CDer(\mathcal{R})$ and $CInn(\mathcal{R})$ are Lie conformal subalgebras of $gc(\mathcal{R})$. Moreover, we have the following tower $CInn(\mathcal{R})\subseteq CDer(\mathcal{R})\subseteq gc(\mathcal{R})$.

\section{construction of conformal triple derivation on Lie conformal superalgebras}
In this section, we discuss the definition and notation of (generalized) conformal triple derivations on Lie conformal superalgebras at first. Then we investigate the conformal algebra structure and some properties of these derivations. For convenience, we always assume that $\mathcal{R}$ is a Lie conformal superalgebra and $\alpha, \beta \in \{\overline{0}, \overline{1}\}$ in this section unless otherwise specified.

\begin{definition}\label{def3.1}
\begin{em}
 A conformal linear map $\phi_{x} \in gc(\mathcal{R})_\alpha$ is called a \emph{conformal triple derivation of degree $\alpha$} of $\mathcal{R}$ if it satisfies the following axiom:
\begin{align}\label{eq3.1}
&\phi_x([[a_{\lambda}b]_{\lambda+\mu}c])\nonumber\\
=&[[\phi_x(a)_{\lambda+x}b]_{\lambda+\mu+x}c]+(-1)^{|\phi||a|}[[a_{\lambda}\phi_x(b)]_{\lambda+\mu+x}c]+(-1)^{|\phi|(|a|+|b|)}[[a_{\lambda}b]_{\lambda+\mu}\phi_x(c)],
\end{align}
for any homogeneous elements $a,b,c \in \mathcal{R}$.

A conformal linear map $\phi_{x} \in gc(\mathcal{R})_\alpha$ is said to be a \emph{generalized conformal triple derivation of degree $\alpha$} of $\mathcal{R}$ if there exists a conformal triple derivation $\tau_x$ of degree $\alpha$ of $\mathcal{R}$ such that
\begin{align}\label{eq3.2}
&\phi_x([[a_{\lambda}b]_{\lambda+\mu}c])\nonumber\\
=&[[\phi_x(a)_{\lambda+x}b]_{\lambda+\mu+x}c]+(-1)^{|\tau||a|}[[a_{\lambda}\tau_x(b)]_{\lambda+\mu+x}c]+(-1)^{|\tau|(|a|+|b|)}[[a_{\lambda}b]_{\lambda+\mu}\tau_x(c)],
\end{align}
for any homogeneous elements $a,b,c \in \mathcal{R}$, and $\tau_x$ is called the \emph{relating conformal triple derivation}, 
\end{em}		
\end{definition}
Denote by $CTDer(\mathcal{R})$ and $GCTDer(\mathcal{R})$ the sets of conformal triple derivations and generalized conformal triple derivations of $\mathcal{R}$ respectively. Obviously, conformal triple derivations are all generalized conformal triple derivations, i.e. $CTDer(\mathcal{R})\subseteq GCTDer(\mathcal{R})$ by setting $\phi_{x}=\tau_x$ in (\ref{eq3.2}). However, the converse is not true in general. 

By above definition, we can obtain the following result immediately.
\begin{lemma}\label{lm3.1}
 For any homogeneous elements $a,b,c \in \mathcal{R}$, $\phi_x \in GCTDer(\mathcal{R})_\alpha$ and $\tau_x$ is the relating conformal triple derivation of $\phi_x$, we can obtain that
\begin{enumerate}
\item $\phi_x([a_{\lambda}[b_{\mu}c]])=[\phi_x(a)_{\lambda+x}[b_{\mu}c]]+(-1)^{\alpha|a|}[a_{\lambda}[\tau_x(b)_{\mu+x}c]]+(-1)^{\alpha(|a|+|b|)}[a_{\lambda}[b_{\mu}\tau_x(c)]].$
\item $(\phi_x-\tau_x)([[a_{\lambda}b]_{\lambda+\mu}c])=[[(\phi_x-\tau_x)(a)_{\lambda+x}b]_{\lambda+\mu+x}c]=(-1)^{\alpha|a|}[[a_{\lambda}(\phi_x-\tau_x)(b)]_{\lambda+\mu+x}c]\\
=(-1)^{\alpha(|a|+|b|)}[[a_{\lambda}b]_{\lambda+\mu}(\phi_x-\tau_x)(c)]$.
\item $(\phi_x-\tau_x)([a_{\lambda}[b_{\mu}c]])=[(\phi_x-\tau_x)(a)_{\lambda+x}[b_{\mu}c]]=(-1)^{\alpha|a|}[a_{\lambda}[(\phi_x-\tau_x)(b)_{\mu+x}c]]\\=(-1)^{\alpha(|a|+|b|)}[a_{\lambda}[b_{\mu}(\phi_x-\tau_x)(c)]].$
\end{enumerate}
\end{lemma}
\begin{proof}
(1) By Jacobi identity, we can obtain that
\begin{align*}
\phi_x([a_{\lambda}[b_{\mu}c]])=&\phi_x([[a_\lambda b]_{\lambda+\mu} c])+(-1)^{|a||b|}\phi_x([ b_\mu [a_\lambda c]])\\
=&\phi_x([[a_\lambda b]_{\lambda+\mu} c])-(-1)^{|a||b|+(|a|+|c|)|b|}\phi_x([[a_\lambda c]_{-\partial-\mu}b])\\
=&[[\phi_x(a)_{\lambda+x}b]_{\lambda+\mu+x}c]+(-1)^{\alpha|a|}[[a_{\lambda}\tau_x(b)]_{\lambda+\mu+x}c]\\&+(-1)^{\alpha(|a|+|b|)}[[a_{\lambda}b]_{\lambda+\mu}\tau_x(c)]-(-1)^{|b||c|}([[\phi_x(a)_{\lambda+x}c]_{-\partial-x-\mu+x}b]+\\&(-1)^{\alpha|a|}[[a_{\lambda}\tau_x(c)]_{-\partial-x-\mu+x}b]+(-1)^{\alpha(|a|+|c|)}[[a_{\lambda}c]_{-\partial-x-\mu}\tau_x(b)])\\
=&[[\phi_x(a)_{\lambda+x}b]_{\lambda+\mu+x}c]+(-1)^{\alpha|a|}[[a_{\lambda}\tau_x(b)]_{\lambda+\mu+x}c]\\&+(-1)^{\alpha(|a|+|b|)}[[a_{\lambda}b]_{\lambda+\mu}\tau_x(c)]
-(-1)^{|b||c|}[[\phi_x(a)_{\lambda+x}c]_{-\partial-\mu}b]\\&-(-1)^{\alpha|a|+|b||c|}[[a_{\lambda}\tau_x(c)]_{-\partial-\mu}b]-(-1)^{\alpha(|a|+|c|)+|b||c|}[[a_{\lambda}c]_{-\partial-x-\mu}\tau_x(b)])\\
=&[[\phi_x(a)_{\lambda+x}b]_{\lambda+\mu+x}c]+ (-1)^{\alpha|a|}[[a_{\lambda}\tau_x(b)]_{\lambda+\mu+x}c]\\&+ (-1)^{\alpha(|a|+|b|)}[[a_{\lambda}b]_{\lambda+\mu}\tau_x(c)]
+(-1)^{(\alpha+|a|)|b|}[b_\mu[\phi_x(a)_{\lambda+x}c]]\\&+ (-1)^{\alpha(|a|+|b|)+|a|||b|}[b_\mu[a_{\lambda}\tau_x(c)]]+
(-1)^{|a||b|}[\tau_x(b)_{\mu+x}[a_{\lambda}c]]\\
=&[\phi_x(a)_{\lambda+x}[b_{\mu}c]]+ (-1)^{\alpha|a|}[a_{\lambda}[\tau_x(b)_{\mu+x}c]]+ (-1)^{\alpha(|a|+|b|)}[a_{\lambda}[b_{\mu}\tau_x(c)]].
\end{align*}
(2) Since $\tau$ is the relating conformal triple derivation of $\phi_x$, we have 
 \begin{align*}
\phi_x([[a_{\lambda}b]_{\lambda+\mu}c])&=[[\phi_x(a)_{\lambda+x}b]_{\lambda+\mu+x}c]+(-1)^{|\tau||a|}[[a_{\lambda}\tau_x(b)]_{\lambda+\mu+x}c]+(-1)^{|\tau|(|a|+|b|)}[[a_{\lambda}b]_{\lambda+\mu}\tau_x(c)]\\
&=[[\phi_x(a)_{\lambda+x}b]_{\lambda+\mu+x}c]+\tau_x([[a_{\lambda}b]_{\lambda+\mu}c])-[[\tau_x(a)_{\lambda+x}b]_{\lambda+\mu+x}c].
 \end{align*}
 That is,
  \begin{align*}
  (\phi_x-\tau_x)([[a_{\lambda}b]_{\lambda+\mu}c])=[[(\phi_x-\tau_x)(a)_{\lambda+x}b]_{\lambda+\mu+x}c].
 \end{align*}
 
 Similarly, by skew-symmetry, we can obtain that 
  \begin{align*}
  (\phi_x-\tau_x)([[a_{\lambda}b]_{\lambda+\mu}c])=(-1)^{\alpha|a|}[[a_{\lambda}(\phi_x-\tau_x)(b)]_{\lambda+\mu+x}c].
 \end{align*}
 
 By skew-symmetry and Jacobi identity, we can deduce that 
 \begin{align*}
 &(\phi_x-\tau_x)([[a_{\lambda}b]_{\lambda+\mu}c])\\&= (\phi_x-\tau_x)([a_\lambda [b_\mu c]])-(-1)^{|a||b|} (\phi_x-\tau_x)([ b_\mu [a_\lambda c]])\\
 &=-(-1)^{|a|(|b|+|c|)}(\phi_x-\tau_x)([[b_\mu c]_{-\partial-\lambda}a])+(-1)^{|a||b|+|b|(|a|+|c|)}(\phi_x-\tau_x)([[a_\lambda c]_{-\partial-\mu}b])\\
 &=-(-1)^{|a|(|b|+|c|)+\alpha|b|}[[b_\mu (\phi_x-\tau_x)(c)]_{-\partial-\lambda}a]+(-1)^{|b||c|+\alpha|a|}[[a_\lambda (\phi_x-\tau_x)(c)]_{-\partial-\mu}b]\\
 &=(-1)^{\alpha(|a|+|b|)}[a_\lambda[b_\mu (\phi_x-\tau_x)(c)]]-(-1)^{\alpha(|a|+|b|)}(-1)^{|a||b|}[b_\mu[a_\lambda (\phi_x-\tau_x)(c)]]\\
 &=(-1)^{\alpha(|a|+|b|)}[[a_{\lambda}b]_{\lambda+\mu}(\phi_x-\tau_x)(c)].
 \end{align*}
 Hence, we have
 \begin{align*}
(\phi_x-\tau_x)([[a_{\lambda}b]_{\lambda+\mu}c])&=[[(\phi_x-\tau_x)(a)_{\lambda+x}b]_{\lambda+\mu+x}c]=(-1)^{\alpha|a|}[[a_{\lambda}(\phi_x-\tau_x)(b)]_{\lambda+\mu+x}c]\\
&=(-1)^{\alpha(|a|+|b|)}[[a_{\lambda}b]_{\lambda+\mu}(\phi_x-\tau_x)(c)].
\end{align*}
(3) It follows immediately from (2) together with the skew-symmetry.

This completes the proof.
\end{proof}

Inspired by the proof of (2) in Lemma \ref{lm3.1}, we can obtain some interesting conclusions as follows:
\begin{corollary}\label{cor3.3}
(1) For any homogeneous elements $a,b,c \in \mathcal{R}$, $\phi_x \in GCTDer(\mathcal{R})_\alpha$, then the following three equations are equivalent:
\begin{enumerate}[(i)]
\item $\phi_x([[a_{\lambda}b]_{\lambda+\mu}c])=[[\phi_x(a)_{\lambda+x}b]_{\lambda+\mu+x}c].$
\item $\phi_x([[a_{\lambda}b]_{\lambda+\mu}c])=(-1)^{\alpha|a|}[[a_{\lambda}\phi_x(b)]_{\lambda+\mu+x}c].$
\item $\phi_x([[a_{\lambda}b]_{\lambda+\mu}c])=[[\phi_x(a)_{\lambda+x}b]_{\lambda+\mu+x}c]=(-1)^{\alpha|a|}[[a_{\lambda}\phi_x(b)]_{\lambda+\mu+x}c]\\=(-1)^{\alpha(|a|+|b|)}[[a_{\lambda}b]_{\lambda+\mu}\phi_x(c)].$
\end{enumerate}
(2) For any homogeneous elements $a,b,c \in \mathcal{R}$, $\phi_x \in GCTDer(\mathcal{R})_\alpha$, then the following three equations are also equivalent:
\begin{enumerate}[(i)]
\item $[[\phi_x(a)_{\lambda+x}b]_{\lambda+\mu+x}c]=(-1)^{\alpha(|a|+|b|)}[[a_{\lambda}b]_{\lambda+\mu}\phi_x(c)].$
\item $[[a_{\lambda}\phi_x(b)]_{\lambda+\mu+x}c]=(-1)^{\alpha|b|}[[a_{\lambda}b]_{\lambda+\mu}\phi_x(c)].$
\item $[[\phi_x(a)_{\lambda+x}b]_{\lambda+\mu+x}c]=(-1)^{\alpha|a|}[[a_{\lambda}\phi_x(b)]_{\lambda+\mu+x}c]=(-1)^{\alpha(|a|+|b|)}[[a_{\lambda}b]_{\lambda+\mu}\phi_x(c)].$
\end{enumerate}
\end{corollary}

\begin{definition}\label{def3.3}
\begin{em}
 A conformal linear map $\phi_{x} \in gc(\mathcal{R})_\alpha$ is called a \emph{conformal triple centroid of degree $\alpha$} of $\mathcal{R}$, if it satisfies the following axiom:
\begin{align}\label{eq3.3}
\phi_x([[a_{\lambda}b]_{\lambda+\mu}c])=[[\phi_x(a)_{\lambda+x}b]_{\lambda+\mu+x}c],
\end{align}
for any homogeneous elements $a,b,c \in \mathcal{R}$.

A conformal linear map $\phi_{x} \in gc(\mathcal{R})_\alpha$ is said to be a \emph{conformal triple quasicentroid of degree $\alpha$} of $\mathcal{R}$, if it satisfies the following equality:
\begin{align}\label{eq3.4}
[[\phi_x(a)_{\lambda+x}b]_{\lambda+\mu+x}c]=(-1)^{\alpha(|a|+|b|)}[[a_{\lambda}b]_{\lambda+\mu}\phi_x(c)],
\end{align}
for any homogeneous elements $a,b,c \in \mathcal{R}$. 

A conformal linear map $\phi_{x} \in gc(\mathcal{R})_\alpha$ is said to be a \emph{central conformal triple derivation of degree $\alpha$} of $\mathcal{R}$, if it satisfies the following equality:
\begin{align}\label{eq3.5}
\phi_x([[a_{\lambda}b]_{\lambda+\mu}c])=[[\phi_x(a)_{\lambda+x}b]_{\lambda+\mu+x}c]=0,
\end{align}
for any homogeneous elements $a,b,c \in \mathcal{R}$. 

\end{em}		
\end{definition}

Throughout what follows, we denote by $TC(\mathcal{R})$, $TQC(\mathcal{R})$ and $ZTDer(\mathcal{R})$ the sets of conformal triple centroids, quasicentroids and central conformal triple derivations respectively. Obviously, we can obtain the tower $ZTDer(\mathcal{R})\subseteq TC(\mathcal{R})\subseteq TQC(\mathcal{R})$. As usual, the center of Lie conformal superalgebra $\mathcal{R}$ is denoted by $Z(\mathcal{R})=\{a\in \mathcal{R} | [a_{\lambda}b]=0, \ \ \forall b\in \mathcal{R} \}$. For a subset $\mathcal{S}$ of $\mathcal{R}$, denote by $C_{\mathcal{R}}(\mathcal{S})$ the centralizer of $\mathcal{S}$ in $\mathcal{R}$. It is clearly that, if $Z(\mathcal{R})=0$, then $ZTDer(\mathcal{R})=0$. 

\begin{remark}
\begin{em}
Inspired by the result of (2) in Lemma \ref{lm3.1}, we can deduce that for any homogeneous elements $a,b,c \in \mathcal{R}$, $\phi_x \in GCTDer(\mathcal{R})_\alpha$ and $\tau_x$ is the relating conformal triple derivation of $\phi_x$, then $\phi_x-\tau_x \in TC(\mathcal{R})_\alpha$.
\end{em} 
\end{remark}

From now on, we can investigate the connections of these conformal triple derivations.
\begin{proposition}\label{pro3.4}
Let $\mathcal{R}$ be a Lie conformal superalgebra. Then $CTDer(\mathcal{R}), GCTDer(\mathcal{R})$ and $TC(\mathcal{R})$ are subalgebras of Lie conformal superalgebra $gc(\mathcal{R})$.
\end{proposition}
\begin{proof}
It is not difficult to see that the proofs of all cases are similar. So we only need to prove that, for any $\phi_x \in CTDer(\mathcal{R})_\alpha$, $\psi_x \in CTDer(\mathcal{R})_\beta$, we have $[\phi_x\psi]_y \in CTDer(\mathcal{R})_{\alpha+\beta}[x]$.

Since $\phi_x \in CTDer(\mathcal{R})_\alpha$, $\psi_x \in CTDer(\mathcal{R})_\beta$, we can obtain that 
\begin{align*}
\phi_x([[a_{\lambda}b]_{\lambda+\mu}c])=[[\phi_x(a)_{\lambda+x}b]_{\lambda+\mu+x}c]+(-1)^{|\phi||a|}[[a_{\lambda}\phi_x(b)]_{\lambda+\mu+x}c]+(-1)^{|\phi|(|a|+|b|)}[[a_{\lambda}b]_{\lambda+\mu}\phi_x(c)],\\
\psi_x([[a_{\lambda}b]_{\lambda+\mu}c])=[[\psi_x(a)_{\lambda+x}b]_{\lambda+\mu+x}c]+(-1)^{|\psi||a|}[[a_{\lambda}\psi_x(b)]_{\lambda+\mu+x}c]+(-1)^{|\psi|(|a|+|b|)}[[a_{\lambda}b]_{\lambda+\mu}\psi_x(c)],
\end{align*}
for any homogeneous elements $a,b,c \in \mathcal{R}$. Furthermore, we can get 
\begin{align*}
\phi_x&\psi_{y-x}([[a_{\lambda}b]_{\lambda+\mu}c])\\
=&\phi_x([[\psi_{y-x}(a)_{\lambda+y-x}b]_{\lambda+\mu+y-x}c]+(-1)^{|\psi||a|}[[a_{\lambda}\psi_{y-x}(b)]_{\lambda+\mu+y-x}c]+ (-1)^{|\psi|(|a|+|b|)}[[a_{\lambda}b]_{\lambda+\mu}\psi_{y-x}(c)])\\
=&[[\phi_x\psi_{y-x}(a)_{\lambda+y}b]_{\lambda+\mu+y}c]+(-1)^{|\phi|(|\psi|+|a|)}[[\psi_{y-x}(a)_{\lambda+y-x}\phi_x(b)]_{\lambda+\mu+y}c]\\&+(-1)^{|\phi|(|\psi|+|a|+|b|)}[[\psi_{y-x}(a)_{\lambda+y-x}b]_{\lambda+\mu+y-x}\phi_x(c)]\\
&+(-1)^{|\psi||a|}([[\phi_x(a)_{\lambda+x}\psi_{y-x}(b)]_{\lambda+\mu+y}c]+(-1)^{|\phi||a|}[[a_{\lambda}\phi_x\psi_{y-x}(b)]_{\lambda+\mu+y}c]\\&+(-1)^{|\phi|(|\psi|+|a|+|b|)}[[a_{\lambda}\psi_{y-x}(b)]_{\lambda+\mu+y-x}\phi_x(c)])\\
&+(-1)^{|\psi|(|a|+|b|)}([[\phi_x(a)_{\lambda+x}b]_{\lambda+\mu+x}\psi_{y-x}(c)]+(-1)^{|\phi||a|}[[a_{\lambda}\phi_x(b)]_{\lambda+\mu+x}\psi_{y-x}(c)]\\&+(-1)^{|\phi|(|a|+|b|)}[[a_{\lambda}b]_{\lambda+\mu}\phi_x\psi_{y-x}(c)])\\
=&[[\phi_x\psi_{y-x}(a)_{\lambda+y}b]_{\lambda+\mu+y}c]+
(-1)^{|\phi|(|\psi|+|a|)}[[\psi_{y-x}(a)_{\lambda+y-x}\phi_x(b)]_{\lambda+\mu+y}c]
\\&+(-1)^{|\phi|(|\psi|+|a|+|b|)}[[\psi_{y-x}(a)_{\lambda+y-x}b]_{\lambda+\mu+y-x}\phi_x(c)]\\
&+(-1)^{|\psi||a|}[[\phi_x(a)_{\lambda+x}\psi_{y-x}(b)]_{\lambda+\mu+y}c]
+(-1)^{(|\phi|+|\psi|)|a|}[[a_{\lambda}\phi_x\psi_{y-x}(b)]_{\lambda+\mu+y}c]
\\&+(-1)^{(|\phi|+|\psi|)|a|+|\phi|(|\psi|+|b|)}[[a_{\lambda}\psi_{y-x}(b)]_{\lambda+\mu+y-x}\phi_x(c)]\\
&+(-1)^{|\psi|(|a|+|b|)}[[\phi_x(a)_{\lambda+x}b]_{\lambda+\mu+x}\psi_{y-x}(c)]
+(-1)^{(|\phi|+|\psi|)|a|+|\psi||b|}[[a_{\lambda}\phi_x(b)]_{\lambda+\mu+x}\psi_{y-x}(c)]\\
&+(-1)^{(|\phi|+|\psi|)(|a|+|b|)}[[a_{\lambda}b]_{\lambda+\mu}\phi_x\psi_{y-x}(c)]\\
=&[[\phi_x\psi_{y-x}(a)_{\lambda+y}b]_{\lambda+\mu+y}c]+(-1)^{(|\phi|+|\psi|)|a|}[[a_{\lambda}\phi_x\psi_{y-x}(b)]_{\lambda+\mu+y}c]\\&+(-1)^{(|\phi|+|\psi|)(|a|+|b|)}[[a_{\lambda}b]_{\lambda+\mu}\phi_x\psi_{y-x}(c)]\\
&+(-1)^{|\phi|(|\psi|+|a|)}[[\psi_{y-x}(a)_{\lambda+y-x}\phi_x(b)]_{\lambda+\mu+y}c]
+(-1)^{|\phi|(|\psi|+|a|+|b|)}[[\psi_{y-x}(a)_{\lambda+y-x}b]_{\lambda+\mu+y-x}\phi_x(c)]\\&+(-1)^{|\psi||a|}[[\phi_x(a)_{\lambda+x}\psi_{y-x}(b)]_{\lambda+\mu+y}c]+(-1)^{(|\phi|+|\psi|)|a|+|\phi|(|\psi|+|b|)}[[a_{\lambda}\psi_{y-x}(b)]_{\lambda+\mu+y-x}\phi_x(c)]\\&+(-1)^{|\psi|(|a|+|b|)}[[\phi_x(a)_{\lambda+x}b]_{\lambda+\mu+x}\psi_{y-x}(c)]
+(-1)^{(|\phi|+|\psi|)|a|+|\psi||b|}[[a_{\lambda}\phi_x(b)]_{\lambda+\mu+x}\psi_{y-x}(c)]
\end{align*}
Similarly, we have
\begin{align*}
\psi_{y-x}&\phi_x([[a_{\lambda}b]_{\lambda+\mu}c])\\
=&[[\psi_{y-x}\phi_x(a)_{\lambda+y}b]_{\lambda+\mu+y}c]+(-1)^{(|\phi|+|\psi|)|a|}[[a_{\lambda}\psi_{y-x}\phi_x(b)]_{\lambda+\mu+y}c]\\&+(-1)^{(|\phi|+|\psi|)(|a|+|b|)}[[a_{\lambda}b]_{\lambda+\mu}\psi_{y-x}\phi_x(c)]\\&
+(-1)^{|\psi|(|\phi|+|a|)}[[\phi_x(a)_{\lambda+x}\psi_{y-x}(b)]_{\lambda+\mu+y}c]
+(-1)^{|\psi|(|\phi|+|a|+|b|)}[[\phi_x(a)_{\lambda+x}b]_{\lambda+\mu+x}\psi_{y-x}(c)]\\
&+(-1)^{|\phi||a|}[[\psi_{y-x}(a)_{\lambda+y-x}\phi_x(b)]_{\lambda+\mu+y}c]
+(-1)^{(|\phi|+|\psi|)|a|+|\psi|(|\phi|+|b|)}[[a_{\lambda}\phi_x(b)]_{\lambda+\mu+x}\psi_{y-x}(c)]\\
&+(-1)^{|\phi|(|a|+|b|)}[[\psi_{y-x}(a)_{\lambda+y-x}b]_{\lambda+\mu+y-x}\phi_x(c)]
+(-1)^{(|\psi|+|\phi|)|a|+|\phi||b|)}[[a_{\lambda}\psi_{y-x}(b)]_{\lambda+ \mu+ y- x}\phi_x(c)].
\end{align*}
Thus, we can obtain that
\begin{align*}
&[\phi_x\psi]_y([[a_{\lambda}b]_{\lambda+\mu}c])=[[[\phi_x\psi]_y(a)_{\lambda+y}b]_{\lambda+\mu+y}c]\\&+(-1)^{(|\phi|+|\psi|)|a|}[[a_{\lambda}[\phi_x\psi]_y(b)]_{\lambda+\mu+y}c]+(-1)^{(|\phi|+|\psi|)(|a|+|b|)}[[a_{\lambda}b]_{\lambda+\mu}[\phi_x\psi]_y(c)].
\end{align*}
Hence, $[\phi_x\psi]_y \in CTDer(\mathcal{R})_{\alpha+\beta}[x]$, i.e. $CTDer(\mathcal{R})$ is a Lie conformal subalgebra of $gc(\mathcal{R})$.

This completes the proof.
\end{proof}

\begin{proposition}\label{pro3.5}
Let $\mathcal{R}$ be a Lie conformal superalgebra. Then $ZTDer(\mathcal{R})$ is a Lie conformal ideal of $CTDer(\mathcal{R})$ and $GCTDer(\mathcal{R})$.
\end{proposition}
\begin{proof}
Obviously, we have the tower $ZTDer(\mathcal{R})\subseteq CTDer(\mathcal{R})\subseteq GCTDer(\mathcal{R})$. We only need to prove that, for any $\phi_x\in ZTDer(\mathcal{R})_\alpha$ and $\psi_x\in CTDer(\mathcal{R})_\beta$, then $[\phi_x\psi]_y \in ZTDer(\mathcal{R})_{\alpha+\beta}[x]$. The other case can be proved similarly.

For any homogeneous elements $a,b,c \in \mathcal{R}$, we have
\begin{align*}
[[[\phi_x\psi]_y&(a)_{\lambda+y}b]_{\lambda+\mu+y}c]\\=&[[\phi_x\psi_{y-x}(a)_{\lambda+y}b]_{\lambda+\mu+y}c]-(-1)^{|\phi||\psi|}[[\psi_{y-x}\phi_x(a)_{\lambda+y}b]_{\lambda+\mu+y}c]\\
=&0-(-1)^{|\phi||\psi|}(\psi_{y-x}([[\phi_x(a)_{\lambda+x}b]_{\lambda+\mu+x}c])\\
&-(-1)^{|\psi|(|\phi|+|a|)}[[\phi_x(a)_{\lambda+x}\psi_{y-x}(b)]_{\lambda+\mu+y}c]-(-1)^{|\psi|(|\phi|+|a|+|b|)}[[\phi_x(a)_{\lambda+x}b]_{\lambda+\mu+x}\psi_{y-x}(c)])\\
=&0.
\end{align*}
Similarly, we can obtain that 
\begin{align*}
[\phi_x\psi]_y([[a_{\lambda}b]_{\lambda+\mu}c])=0.
\end{align*}
Hence, we can deduce that $[\phi_x\psi]_y \in ZTDer(\mathcal{R})_{\alpha+\beta}[x]$.

This completes the proof.
\end{proof}

\begin{remark}
\begin{em}
It is well known that $gc(\mathcal{R})$ is a simple Lie conformal superalgebra. Thus, if $GCTDer(\mathcal{R})=gc(\mathcal{R})$, by Proposition \ref{pro3.5}, $ZTDer(\mathcal{R})$ is equal to either 0 or $gc(\mathcal{R})$. Moreover, if $ZTDer(\mathcal{R})=gc(\mathcal{R})$, then $\mathcal{R}$ is abelian.
\end{em}
\end{remark}

\begin{proposition}\label{pro3.7}
Let $\mathcal{R}$ be a Lie conformal superalgebra. If $Z(\mathcal{R})=0$, then $TC(\mathcal{R})$ and $TQC(\mathcal{R})$ are commutative Lie conformal algebras. In particular, if $Z(\mathcal{R})=0$, $TC(\mathcal{R})$ centralizes $TQC(\mathcal{R})$.
\end{proposition}
\begin{proof}
Clearly, we only need to prove the first case, the other case can be proved similarly. For any homogeneous elements $a,b,c \in \mathcal{R}$, $\phi_x\in TC(\mathcal{R})_\alpha$, $\psi_x\in TC(\mathcal{R})_\beta$, we can deduce that 
\begin{align*}
[[\phi_x\psi_{y-x}(a)_{\lambda+y}b]_{\lambda+\mu+y}c]&=\phi_{x}([[\psi_{y-x}(a)_{\lambda+y-x}b]_{\lambda+\mu+y-x}c])\\
&=(-1)^{|\psi||a|}\phi_x([[a_{\lambda}\psi_{y-x}(b)]_{\lambda+\mu+y-x}c])\\
&=(-1)^{|\psi||a|}[[\phi_x(a)_{\lambda+x}\psi_{y-x}(b)]_{\lambda+\mu+y}c]\\
&=(-1)^{|\psi||a|+|\psi|(|\phi|+|a|)}[[\psi_{y-x}\phi_x(a)_{\lambda+y}b]_{\lambda+\mu+y}c]\\
&=(-1)^{|\psi||\phi|}[[\psi_{y-x}\phi_x(a)_{\lambda+y}b]_{\lambda+\mu+y}c]\\.
\end{align*}
Hence, we have 
\begin{align*}
[[[\phi_x\psi]_y(a)_{\lambda+y}b]_{\lambda+\mu+y}c]=0.
\end{align*}
As $Z(\mathcal{R})=0$ and  $b,c$ are arbitrary elements of $\mathcal{R}$, so we can find that
\begin{align*}
[\phi_x\psi]_y(a)=0.
\end{align*}
Thus, we have $[\phi_x\psi]_y=0$.

This completes the proof.
\end{proof}

\begin{proposition}\label{pro3.8}
Let $\mathcal{R}$ be a Lie conformal algebra. Then we have the following results:
\begin{enumerate}
\item $[{CTDer(\mathcal{R})_{\alpha}}_\lambda TC(\mathcal{R})_{\beta}]_\mu \subseteq TC(\mathcal{R})_{\alpha+\beta}[\lambda]$.
\item $[{GCTDer(\mathcal{R})_{\alpha}}_\lambda TC(\mathcal{R})_{\beta}]_\mu \subseteq TC(\mathcal{R})_{\alpha+\beta}[\lambda]$.
\item $[{CTDer(\mathcal{R})_{\alpha}}_\lambda TQC(\mathcal{R})_{\beta}]_\mu \subseteq TQC(\mathcal{R})_{\alpha+\beta}[\lambda]$.
\end{enumerate}
\end{proposition}
\begin{proof}

Obviously, we only need to prove the first case, other cases can be proved similarly.
Since $\phi_x\in CTDer(\mathcal{R})_{\alpha}$, $\psi_x \in TC(\mathcal{R})_{\beta}$, for any homogeneous elements $a,b,c \in \mathcal{R}$, we can obtain that 
\begin{align*}
&[[\phi_x\psi_{y-x}(a)_{\lambda+y}b]_{\lambda+\mu+y}c]\\
=&\phi_x([[\psi_{y-x}(a)_{\lambda+y-x}b]_{\lambda+\mu+y-x}c])-(-1)^{|\phi|(|\psi|+|a|)}[[\psi_{y-x}(a)_{\lambda+y-x}\phi_x(b)]_{\lambda+\mu+y}c]\\
&-(-1)^{|\phi|(|\psi|+|a|+|b|)}[[\psi_{y-x}(a)_{\lambda+y-x}b]_{\lambda+\mu+y-x}\phi_x(c)]\\
=&\phi_x\psi_{y-x}([[a_{\lambda}b]_{\lambda+\mu}c])-(-1)^{|\phi|(|\psi|+|a|)}\psi_{y-x}([[a_{\lambda}\phi_x(b)]_{\lambda+\mu+x}c])\\&-(-1)^{|\phi|(|\psi|+|a|+|b|)}\psi_{y-x}([[a_{\lambda}b]_{\lambda+\mu}\phi_x(c)]).
\end{align*}
Similarly, we have
\begin{align*}
&[[\psi_{y-x}\phi_x(a)_{\lambda+y}b]_{\lambda+\mu+y}c]\\
=&\psi_{y-x}([[\phi_x(a)_{\lambda+x}b]_{\lambda+\mu+x}c])\\
=&\psi_{y-x}\phi_x([[a_{\lambda}b]_{\lambda+\mu}c])-(-1)^{|\phi||a|}\psi_{y-x}([[a_{\lambda}\phi_x(b)]_{\lambda+\mu+x}c])\\&-(-1)^{|\phi|(|a|+|b|)}\psi_{y-x}([[a_{\lambda}b]_{\lambda+\mu}\phi_x(c)]).
\end{align*}
Hence, we can obtain that
\begin{align*}
[[[\phi_x\psi]_y(a)_{\lambda+y}b]_{\lambda+\mu+y}c]=[\phi_x\psi]_y([[a_{\lambda}b]_{\lambda+\mu}c]).
\end{align*}
Thus, we have $[\phi_x\psi]_y \in TC(\mathcal{R})_{\alpha+\beta}[x]$.

This completes the proof.
\end{proof}
\section{Classification of conformal triple derivations on simple Lie conformal superalgebras}
In this section, we classify (generalized) conformal triple derivations of all simple Lie conformal superalgebras. For convenience, we always assume that $\mathcal{R}$ is a finite simple Lie conformal superalgebra and $\alpha, \beta \in \{\overline{0}, \overline{1}\}$.

It is easy to see that $CDer(\mathcal{R})$ and $CInn(\mathcal{R})$ are Lie conformal subalgebras of $CTDer(\mathcal{R})$. Furthermore, we have the following lemma.
\begin{lemma}\label{lm4.1}
Let $\mathcal{R}$ be a finite simple Lie conformal superalgebra. Then $CInn(\mathcal{R})$ is an ideal of Lie conformal superalgebra $CTDer(\mathcal{R})$.
\end{lemma}
\begin{proof}
Let $\phi_x \in CTDer(\mathcal{R})_\alpha$, $a \in \mathcal{R}_\beta$. Since $\mathcal{R}$ is a simple Lie conformal superalgebra, then there exists some finite index set $I, J \subseteq \mathbb{Z^+}$ and homogeneous elements $a_{ij}^1, a_{ij}^2 \in \mathcal{R}$ such that $a= \sum_{i\in I, j\in J}\frac{d^i}{d\lambda^i}[{a_{ij}^1}_\lambda {a_{ij}^2}]|_{\lambda=0}$.

For any $b\in \mathcal{R}$, we have 
\begin{align*}
[\phi_x ad a]_y(b)&=\phi_x ada_{y-x}(b)-(-1)^{|\phi||a|}ada_{y-x}\phi_x(b)\\
=&\phi_x ([a_{y-x}b])-(-1)^{|\phi||a|}[a_{y-x}\phi_x(b)]\\
=&\phi_x([(\sum_{i\in I, j\in J}\frac{d^i}{d\lambda^i}[{a_{ij}^1}_\lambda {a_{ij}^2}]|_{\lambda=0})_{y-x}b])-(-1)^{|\phi||a|}[(\sum_{i\in I, j\in J}\frac{d^i}{d\lambda^i}[{a_{ij}^1}_\lambda {a_{ij}^2}]|_{\lambda=0})_{y-x}\phi_x(b)]\\
=&\sum_{i\in I, j\in J}\frac{d^i}{d\lambda^i}([[\phi_x(a_{ij}^1)_{\lambda+ x} {a_{ij}^2}]_{y}b]+(-1)^{|\phi||a_{ij}^1|}[[{a_{ij}^1}_\lambda \phi_x({a_{ij}^2})]_{y}b]\\+&(-1)^{|\phi|(|a_{ij}^1|+|a_{ij}^2|)}[[{a_{ij}^1}_\lambda {a_{ij}^2}]_{y-x}\phi_x(b)])|_{\lambda=0}
-(-1)^{|\phi||a|}[(\sum_{i\in I, j\in J}{\frac{d^i}{d\lambda^i}[{a_{ij}^1}_\lambda {a_{ij}^2}]|_{\lambda=0}})_{y-x}\phi_x(b)]\\
=&\sum_{i\in I, j\in J}\frac{d^i}{d\lambda^i}([[\phi_x(a_{ij}^1)_{\lambda+ x} {a_{ij}^2}]_{y}b]+(-1)^{|\phi||a_{ij}^1|}[[{a_{ij}^1}_\lambda \phi_x({a_{ij}^2})]_{y}b])|_{\lambda=0}\\
=&ad(\sum_{i\in I, j\in J}\frac{d^i}{d\lambda^i}([\phi_x(a_{ij}^1)_{\lambda+ x} {a_{ij}^2}]+(-1)^{|\phi||a_{ij}^1|}[{a_{ij}^1}_\lambda \phi_x({a_{ij}^2})])|_{\lambda=0})_{y}(b).
\end{align*}
By the arbitrariness of $b$, $[\phi_x ad a]_y$ is an inner conformal derivation. Thus, $CInn(\mathcal{R})$ is an ideal of Lie conformal algebra $CTDer(\mathcal{R})$.
\end{proof}
In fact, there are other connections between $CDer(\mathcal{R})$ and $CTDer(\mathcal{R})$, which is stronger than the relationship between them as shown in above lemma. Thus, we have the following lemma.

\begin{lemma}\label{lm4.2}
Let $\mathcal{R}$ be a finite simple Lie conformal superalgebra. Then there exists a $\mathbb{C}$-linear map $\delta: CTDer(\mathcal{R})\to CDer(\mathcal{R})$, $\phi_x \mapsto \delta_{\phi_x}$ such that for all $a \in \mathcal{R}$, $\phi_x \in CTDer(\mathcal{R})$, one has $[\phi_x ad a]_y=ad {\delta_{\phi_x}}(a)_y$.
\end{lemma}
\begin{proof}
By the proof of Lemma \ref{lm4.1}, take $\phi_x \in CTDer(\mathcal{R})_\alpha$, then we can define a conformal linear endomorphism $\delta_{\phi_x}$ on $\mathcal{R}$, such that for any homogeneous element $a= \sum_{i\in I, j\in J}\frac{d^i}{d\lambda^i}[{a_{ij}^1}_\lambda {a_{ij}^2}]|_{\lambda=0} \in \mathcal{R}$, 
\begin{align*}
\delta_{\phi_x}(a)= \sum_{i\in I, j\in J}\frac{d^i}{d\lambda^i}([\phi_x(a_{ij}^1)_{\lambda+ x} {a_{ij}^2}]+ (-1)^{|\phi||a_{ij}^1|}[ {a_{ij}^1}_\lambda \phi_x({a_{ij}^2})])|_{\lambda=0}.
\end{align*}
At first, $\delta_{\phi_x}(a)$ does not depend on the choice of expression of $a$. For proving it, take 
\begin{align*}
A&= \sum_{i\in I, j\in J}\frac{d^i}{d\lambda^i}([\phi_x(a_{ij}^1)_{\lambda+ x} {a_{ij}^2}]+ (-1)^{|\phi||a_{ij}^1|}[ {a_{ij}^1}_\lambda \phi_x({a_{ij}^2})])|_{\lambda=0},\\
B&= \sum_{s\in I^{'}, t\in J^{'}}\frac{d^s}{d\lambda^s}([\phi_x(b_{st}^1)_{\lambda+ x} {b_{st}^2}]+ (-1)^{|\phi||b_{st}^1|}[ {b_{st}^1}_\lambda \phi_x({b_{st}^2})])|_{\lambda=0},
\end{align*}
where $a$ can also be expressed in the form $a=  \sum_{s\in I^{'}, t\in J^{'}}\frac{d^s}{d\lambda^s}[{b_{st}^1}_\lambda {b_{st}^2}]|_{\lambda=0}$. Since $\phi_x \in CTDer(\mathcal{R})_\alpha$, for any $c \in \mathcal{R}$, we can deduce that 
\begin{align*}
[A_{\lambda+\mu+x}c]=\phi_x([a_{\lambda+\mu}c])-(-1)^{|\phi||a|}[a_{\lambda+\mu}\phi_x(c)]=[B_{\lambda+\mu+x}c].
\end{align*}
Hence, $[(A-B)_{\lambda+\mu+x}c]=0$, for any $c \in \mathcal{R}$, i.e. $A-B\in Z(\mathcal{R})$. Obviously, we have $Z(\mathcal{R})=0$, because $\mathcal{R}$ is a simple Lie conformal superalgebra. Then, we can obtain that $A= B$. Therefore, $\delta_{\phi_x}$ is well-defined. Furthermore, it follows from the proof of Lemma \ref{lm4.1} immediately that $[\phi_x ad a]_y=ad {\delta_{\phi_x}}(a)_y$.

Next, we need to prove that $\delta_{\phi_x}$ is a conformal linear map, i.e. $\delta_{\phi_x}(\partial a)=(\partial+x)\delta_{\phi_x}(a)$. For any homogeneous element $a= \sum_{i\in I, j\in J}\frac{d^i}{d\lambda^i}[{a_{ij}^1}_\lambda {a_{ij}^2}]|_{\lambda=0} \in \mathcal{R}$, we have $\partial a= \sum_{i\in I, j\in J}\frac{d^i}{d\lambda^i}([{\partial a_{ij}^1}_\lambda {a_{ij}^2}] + [{a_{ij}^1}_\lambda \partial {a_{ij}^2}])|_{\lambda=0}$. Thus, we can obtain that
\begin{align*}
\delta_{\phi_x}(\partial a)
	=& \sum_{i\in I, j\in J}\frac{d^i}{d\lambda^i}([\phi_x(\partial a_{ij}^1)_{\lambda+ x} {a_{ij}^2}]+ (-1)^{|\phi||a_{ij}^1|}[{\partial a_{ij}^1}_\lambda \phi_x({a_{ij}^2})]\\&+ [\phi_x(a_{ij}^1)_{\lambda+ x} {\partial  a_{ij}^2}]+ (-1)^{|\phi||a_{ij}^1|}[{a_{ij}^1}_\lambda \phi_x({\partial  a_{ij}^2})])|_{\lambda=0}\\
	=& \sum_{i\in I, j\in J}\frac{d^i}{d\lambda^i}(- \lambda([\phi_x(a_{ij}^1)_{\lambda+ x} {a_{ij}^2}]+ (-1)^{|\phi||a_{ij}^1|}[{a_{ij}^1}_\lambda \phi_x({a_{ij}^2})])\\
	&+ (\partial +\lambda+ x)([\phi_x(a_{ij}^1)_{\lambda+ x} { a_{ij}^2}]+ (-1)^{|\phi||a_{ij}^1|}[{a_{ij}^1}_\lambda \phi_x({a_{ij}^2})]))|_{\lambda=0}\\
	=& (\partial+ x)\sum_{i\in I, j\in J}\frac{d^i}{d\lambda^i}([\phi_x(a_{ij}^1)_{\lambda+ x} { a_{ij}^2}]+ (-1)^{|\phi||a_{ij}^1|}[{a_{ij}^1}_\lambda \phi_x({a_{ij}^2})])|_{\lambda=0}\\
	=& (\partial+ x)\delta_{\phi_x}(a).
\end{align*}

Finally, we should prove that $\delta_{\phi_x}$ is a conformal derivation. By the above discussion, for any  $\phi_x \in CTDer(\mathcal{R})_\alpha$, homogeneous elements $a,b \in \mathcal{R}$, we have $[\phi_x ad[a_\lambda b]]=ad \delta_{\phi_x} ([a_\lambda b])$. On the other hand, we have
\begin{align*}
[\phi_x ad[a_\lambda b]]&=[\phi_x [ad a_\lambda ad b]]\\
&=[[\phi_x ad a]_{\lambda+x} ad b]+(-1)^{|\phi||a|}[ad a_\lambda [\phi_x ad b]]\\
&=[ad \delta_{\phi_x}(a)_{\lambda+x} ad b]+(-1)^{|\phi||a|}[ad a_{\lambda} ad \delta_{\phi_x}(b)]\\
&=ad ([\delta_{\phi_x}(a)_{\lambda+x} b]+(-1)^{|\phi||a|}[a_{\lambda} \delta_{\phi_x}(b)]).
\end{align*}
Therefore, $ad \delta_{\phi_x} ([a_\lambda b])=ad ([\delta_{\phi_x}(a)_{\lambda+x} b]+(-1)^{|\phi||a|}[a_{\lambda} \delta_{\phi_x}(b)])$. Since $Z(\mathcal{R})=0$, then we have $\delta_{\phi_x} ([a_\lambda b])=[\delta_{\phi_x}(a)_{\lambda+x} b]+(-1)^{|\phi||a|}[a_{\lambda} \delta_{\phi_x}(b)]$. By the arbitrariness of $a,b$, we can deduce that $\delta_{\phi_x}  \in  CDer(\mathcal{R})$.

This completes the proof.
\end{proof}

\begin{lemma}\label{lm4.3}
Let $\mathcal{R}$ be a finite simple Lie conformal superalgebra. Then the centralizer of $CInn(\mathcal{R})$ in $CTDer(\mathcal{R})$ is trivial, i.e., $C_{CTDer(\mathcal{R})}(CInn(\mathcal{R}))=0$. In particular, the center of $CTDer(\mathcal{R})$ is zero.
\end{lemma}
\begin{proof}
Suppose that $\phi_x \in C_{CTDer(\mathcal{R})}(CInn(\mathcal{R}))_\alpha$. Then we have $[\phi_x ad a]=0$ for any homogeneous element $a \in \mathcal{R}$. Thus, for any homogeneous element $b \in \mathcal{R}$, we can obtain that $\phi_x ([a_{\lambda}b])-(-1)^{|\phi||a|}[a_\lambda \phi_x (b)]=[\phi_x ad a]_{\lambda+x}(b)=0$. Moreover, by the skew-symmetry we have $\phi_x ([a_{\lambda}b])=(-1)^{|\phi||a|}[a_\lambda \phi_x (b)]=[\phi_x (a)_{\lambda+x} b]$.

On the one hand, we have
\begin{align*}
\phi_x([[a_{\lambda}b]_{\lambda+\mu}c])=[[\phi_x(a)_{\lambda+x}b]_{\lambda+\mu+x}c]=(-1)^{|\phi||a|}[[a_{\lambda}\phi_x(b)]_{\lambda+\mu+x}c]=(-1)^{|\phi|(|a|+|b|)}[[a_{\lambda}b]_{\lambda+\mu}\phi_x(c)],
\end{align*}
for any homogeneous elements $a,b,c \in \mathcal{R}$. On the other hand, we have
\begin{align*}
\phi_x([[a_{\lambda}b]_{\lambda+\mu}c])=[[\phi_x(a)_{\lambda+x}b]_{\lambda+\mu+x}c]+(-1)^{|\phi||a|}[[a_{\lambda}\phi_x(b)]_{\lambda+\mu+x}c]+(-1)^{|\phi|(|a|+|b|)}[[a_{\lambda}b]_{\lambda+\mu}\phi_x(c)],
\end{align*}
because $\phi_x \in CTDer(\mathcal{R})$. Hence, we can obtain that
\begin{align*}
\phi_x([[a_{\lambda}b]_{\lambda+\mu}c])=3\phi_x([[a_{\lambda}b]_{\lambda+\mu}c]). 
\end{align*}
This means $\phi_x([[a_{\lambda}b]_{\lambda+\mu}c])=0$ for any homogeneous elements $a,b,c \in \mathcal{R}$. Since $\mathcal{R}$ is a simple Lie conformal superalgebra, every element of $\mathcal{R}$ can be expressed as the linear combination of finite elements of the form $\sum_{i\in I, j\in J, k \in K}\frac{d^i}{d\lambda^i}\frac{d^k}{d\mu^k}[[{a_{ijk}^1}_\lambda {a_{ijk}^2}]_\mu {a_{ijk}^3} ]|_{\lambda=\mu=0}$, then we can deduce that $\phi_x=0$.
\end{proof}

Now, let us introduce the main result in this section by the following theorem, which will give the classification of the conformal triple derivation of a finite simple Lie conformal superalgebra. 
\begin{theorem}\label{tm4.4}
Let $\mathcal{R}$ be a finite simple Lie conformal superalgebra. Then $CTDer(\mathcal{R})=CDer(\mathcal{R})$.
\end{theorem}
\begin{proof}
Assume that $\phi_x \in CTDer(\mathcal{R})_\alpha$ and $a$ is homogeneous element in $\mathcal{R}$. By Lemma \ref{lm4.2}, we have $[\phi_x ad a]_y=ad {\delta_{\phi_x}}(a)_y$ and $\delta_{\phi_x} \in CDer(\mathcal{R})$. Moreover, for any $b \in \mathcal{R}$, we have $ad {\delta_{\phi_x}}(a)_y (b)=[\delta_{\phi_x}(a)_y b]=\delta_{\phi_x}([a_{y-x}b])-(-1)^{|\phi||a|}[a_{y-x}\delta_{\phi_x}(b)]=[\delta_{\phi_x} ada]_y(b)$. By the arbitrariness of $b$, we can deduce that $ad {\delta_{\phi_x}}(a)_y=[\delta_{\phi_x} ada]_y$. Hence, we have $[\phi_x ad a]_y=[\delta_{\phi_x} ada]_y$, i.e., $\phi_x-\delta_{\phi_x} \in C_{CTDer(\mathcal{R})}(CInn(\mathcal{R}))$. By Lemma \ref{lm4.3}, $\phi_x-\delta_{\phi_x}=0$, i.e., $\phi_x=\delta_{\phi_x}$. Thus, the theorem follows from Lemma \ref{lm4.2}.

\end{proof}

\section{$(\varPhi, \varPsi)$-triple derivations and (A,B,C,D)-triple derivation of Lie  conformal superalgebras}Our aim in this section is to show some fundamental properties of conformal $(\varPhi, \varPsi)$-triple derivations. Later in this section  we classify (generalized) $(A,B,C,D)$-triple derivations of all simple Lie conformal superalgebras $\mathcal{R}$.
\begin{definition}
Let $\mathcal{R}$ be a Lie conformal superalgebra and $G$ be a subgroup of $Aut(\mathcal{R})$. Then $\phi_{x} \in gc(\mathcal{R})_\alpha$ is called a conformal $(\varPhi, \varPsi)$-triple derivation of degree $\alpha$ of $\mathcal{R}$, if there exists two automorphisms $\varPhi, \varPsi \in G$ such that
\begin{align}
	\phi_{x}([[a_{\lambda}b]_{\lambda+\mu} c]) =& [[\phi_{x}(a)_{\lambda+x}\varPhi(b)]_{\lambda+x+\mu}\varPsi(c)] + (-1)^{|\phi||a|} [[\varPhi(a)_{\lambda}\phi_{x}(b)]_{\lambda+ x+ \mu} \varPsi(c)]\nonumber\\
	&+(-1)^{|\phi| (|a|+|b|)}[[\varPhi(a)_{\lambda}\varPsi(b)]_{\lambda+\mu} \phi_{x}(c)]
	\end{align} 
	for any homogeneous elements $a,b,c \in \mathcal{R}.$
\end{definition}
In this case, $\varPhi$ and $\varPsi$ are called the associated automorphisms of $\phi_{x}$. The set of all conformal $(\varPhi, \varPsi)$-triple derivations of $\mathcal{R}$ is denoted  by $CTDer_{\varPhi, \varPsi}(\mathcal{R})$. It is not hard to see that for $\varPhi = \varPsi = id,$ we have $CTDer_{\varPhi, \varPsi}(\mathcal{R}) = CTDer(\mathcal R)$. For convenience, we will abbreviate $CTDer_{\varPhi, id}(\mathcal{R})$ as $CTDer_{\varPhi}(\mathcal{R})$ in this section.

\begin{proposition}
Let $\mathcal{R}$ be a Lie conformal superalgebra and $G$ be a subgroup of $Aut(\mathcal{R})$. For any automorphisms $\varPhi, \varPsi \in G$, $CTDer_{\varPhi, \varPsi}(\mathcal{R})$ is isomorphic to $CTDer_{\varPsi^{-1}\varPhi}(\mathcal{R})$ as a $C[\partial]$-module. Moreover, 	$rank(CTDer_{\varPhi, \varPsi}(\mathcal{R})) = rank(CTDer_{\varPsi^{-1}\varPhi}(\mathcal{R})).$
\end{proposition}
\begin{proof}Consider a map $\sigma_{\varPsi} : CTDer_{\varPhi, \varPsi}(\mathcal{R}) \to CTDer_{\varPsi^{-1}\varPhi}(\mathcal{R})$ by $\sigma_{\varPsi}(\phi_{x}) = \varPsi^{-1}\phi_{x}$, for any $\phi_{x} \in CTDer_{\varPhi, \varPsi}(\mathcal{R})_\alpha$.\\
For any homogeneous elements $ a,b,c \in \mathcal\mathcal{R}$, we have
\nonumber\begin{align}
\varPsi^{-1}&\phi_{x}([[a_{\lambda}b]_{\lambda+\mu}c])\\=&\varPsi^{-1}([[\phi_{x}(a)_{\lambda+x}\varPhi(b)]_{\lambda+x+\mu}\varPsi(c)]) + (-1)^{|\phi|| a|}\varPsi^{-1}([[\varPhi(a)_{\lambda}\phi_{x}(b)]_{\lambda+x+\mu} \varPsi(c)])\\+&(-1)^{|\phi| (|a|+|b|)}\varPsi^{-1}([[\varPhi(a)_{\lambda}\varPsi(b)]_{\lambda+\mu} \phi_{x}(c)])\\=&[[\varPsi^{-1}\phi_{x}(a)_{\lambda+x}\varPsi^{-1}\varPhi(b)]_{\lambda+x+\mu}c] + (-1)^{|\phi||a|}[[\varPsi^{-1}\varPhi(a)_{\lambda}\varPsi^{-1}\phi_{x}(b)]_{\lambda+x+\mu}c]\\+&(-1)^{|\phi| (|a|+|b|)}[[\varPsi^{-1}\varPhi(a)_{\lambda}(b)]_{\lambda+\mu} \varPsi^{-1}\phi_{x}(c)]
\end{align}
Thus, $\varPsi^{-1}\phi_{x}\in CTDer_{\varPsi^{-1}\varPhi,id}(\mathcal{R})=CTDer_{\varPsi^{-1}\varPhi}(\mathcal{R})$. Then $\sigma_{\varPsi}$ is a well defined map.

Moreover, it is easy to see that
\begin{equation*}
\sigma_{\varPsi}(\phi_{x} + \psi_{x}) = \varPsi^{-1}(\phi_{x} + \psi_{x}) = \varPsi^{-1}(\phi_{x})+\varPsi^{-1}(\psi_{x})  = \sigma_{\varPsi}(\phi_{x}) + \sigma_{\varPsi}(\psi_{x})
\end{equation*}
holds for any $\phi_{x}, \psi_{x}\in CTDer_{\varPhi, \varPsi}(\mathcal{R})$. For any $\phi_{x}\in CTDer_{\varPhi, \varPsi}(\mathcal{R})$ and $a\in\mathcal{R}$, we have 
\begin{equation*}
(\sigma_{\varPsi}(\partial\phi))_{x}a =  \varPsi^{-1}(\partial\phi)_{x}a =  \varPsi^{-1}(-x\phi_{x}(a)) =  -x\varPsi^{-1}(\phi_{x}(a)) = (\partial\sigma_{\varPsi}(\phi))_xa.\\
\end{equation*}
Due to the arbitrariness of  $\phi_{x}$ and $a$, we can deduce that $\sigma_{\varPsi}$ is a $C[\partial]$-module homomorphism. 

Similarly, we can define $\tau_{\varPsi} : CTDer_{\varPsi^{-1}\varPhi}(\mathcal{R}) \to CTDer_{\varPhi, \varPsi}(\mathcal{R}) $ by $\tau_{\varPsi}(\varphi_{x}) = \varPsi\varphi_{x}$.
It is easy to see $\tau_{\varPsi}$ is also a well defined  $C[\partial]$-module homomorphism, it follows that $\sigma_{\varPsi}= \tau_{\varPsi}^{-1}$. Thus $CTDer_{\varPsi^{-1}\varPhi}(\mathcal{R})$ is isomorphic to $CTDer_{\varPsi, \varPhi}(\mathcal{R})$ as a $C[\partial]$-module.
\end{proof}

Inspired by Proposition \ref{pro3.4}, we want to define the structure of conformal superalgebra in $CTDer_{\varPhi, \varPsi}(\mathcal{R})$. Similar to the Proposition 2.2 in \cite{Feng-Zhao-Chen}, it is not difficult to see that we can obtain the following result.
\begin{proposition}
Let $\mathcal{R}$ be a Lie conformal superalgebra and $G$ be a subgroup of $Aut(\mathcal{R})$. For any $\varPhi\in G$, $CTDer_{\varPhi,\varPhi}(\mathcal{R})$ can be viewed as a Lie conformal superalgebra and $CTDer_{\varPhi,\varPhi}(\mathcal{R})$ is isomorphic to  $CTDer(\mathcal{R})$ as Lie conformal superalgebra.
\end{proposition}

\begin{proposition}
Let $\mathcal{R}$ be a Lie conformal superalgebra and $G$ be a subgroup of $Aut(\mathcal{R})$. Suppose that $\varPhi,\varPsi$ are two elements in $G$ such that $(\varPhi-\varPsi)(\mathcal{R}) \in Z(\mathcal{R})$, then
	$CTDer_{\varPhi}(\mathcal{R}) = CTDer_{\varPsi}(\mathcal{R})$. In particular, if $(\varPhi-id_\mathcal{R})(\mathcal{R}) \subseteq Z(\mathcal{R})$, then $CTDer_{\varPhi}(\mathcal{R}) = CTDer(\mathcal{R})$.
\end{proposition}

\begin{proof}
For any $\phi_{x} \in CTDer_{\varPhi}(\mathcal{R})_\alpha$ and homogeneous elements $a,b,c\in \mathcal{R}$ and since $(\varPhi-\varPsi)(\mathcal{R}) \in Z(\mathcal{R})$, we can obtain that
	\begin{align}\nonumber
	[[\phi_{x}(a)_{\lambda+x}(\varPhi-\varPsi)(b)]_{\lambda+\mu+x}c]=[[(\varPhi-\varPsi)(a)_{\lambda}\phi_x(b)]_{\lambda+\mu+x}c]=[[(\varPhi-\varPsi)(a)_{\lambda}b]_{\lambda+\mu}\phi_x(c)]=0,
	\end{align}
	i.e,\begin{align*}
	[[\phi_{x}(a)_{\lambda+x}\varPhi(b)]_{\lambda+\mu+x}c]&=[[\phi_{x}(a)_{\lambda+x}\varPsi(b)]_{\lambda+\mu+x}c],\\
	\nonumber[[\varPhi(a)_{\lambda}\phi_x(b)]_{\lambda+\mu+x}c]&=[[\varPsi(a)_{\lambda}\phi_x(b)]_{\lambda+\mu+x}c],\\
	\nonumber[[\varPhi(a)_{\lambda}b]_{\lambda+\mu}\phi_x(c)]&=[[\varPsi(a)_{\lambda}b]_{\lambda+\mu}\phi_x(c)].
	\end{align*}
Thus  we get
\begin{align*}
	&\phi_{x}([[a_{\lambda}b]_{\lambda+\mu}c])\\
	&= [[\phi_{x}(a)_{\lambda+x}\varPhi(b)]_{\lambda+\mu+x}c]+
	(-1)^{|\phi||a|}[[\varPhi(a)_{\lambda}\phi_x(b)]_{\lambda+\mu+x}c]+(-1)^{|\phi|(|a|+|b|)}
	[[\varPhi(a)_{\lambda}b]_{\lambda+\mu}\phi_x(c)]\\
	&= [[\phi_{x}(a)_{\lambda+x}\varPsi(b)]_{\lambda+\mu+x}c]+
	(-1)^{|\phi||a|}[[\varPsi(a)_{\lambda}\phi_x(b)]_{\lambda+\mu+x}c]+(-1)^{|\phi|(|a|+|b|)}
	[[\varPsi(a)_{\lambda}b]_{\lambda+\mu}\phi_x(c)]
	\end{align*}
	It follows that  $\phi_{x} \in CTDer_{\varPsi}(\mathcal{R})$, which implies $CTDer_{\varPhi}(\mathcal{R})\subseteq CTDer_{\varPsi}(\mathcal{R})$. Similarly, we have $CTDer_{\varPsi}(\mathcal{R})\subseteq CTDer_{\varPhi}(\mathcal{R})$. Hence, $CTDer_{\varPsi}(\mathcal{R})= CTDer_{\varPhi}(\mathcal{R})$.
	Proof of the second part can be done similarly by taking $\varPsi= id_\mathcal{R}$.
\end{proof}

\begin{proposition}\label{prop5.5}
Let $\mathcal{R}$ be a Lie conformal superalgebra and $G$ be a subgroup of $Aut(\mathcal{R})$. Suppose that $\varPhi, \varPsi$  are two elements in $G$ such that $\phi_{x} \in CTDer_{\varPhi}(\mathcal{R})_\alpha$ and $\psi_{x} \in CTDer_{\varPsi} (\mathcal{R})_\beta$. If $\varPhi$ and $\varPsi$ commute, $\phi_{x}$ commutes with $\varPsi$ and $\psi_y$  commutes with $\varPhi$, then $[\phi_{x}\psi]_{y}\in CTDer_{\varPhi\varPsi}(\mathcal{R})_{\alpha+\beta}.$
\end{proposition}

\begin{proof}
For any homogeneous elements $a,b,c \in \mathcal{R}$, since $\varPhi$ and $\varPsi$ commute, $\phi_{x}$ commutes with $\varPsi$ and $\psi_y$  commutes with $\varPhi$, we have 
	\begin{align*}
	\phi_x&\psi_{y- x}([[a_{\lambda}b]_{\lambda+ \mu}c])\\
	=& \phi_x([[\psi_{y- x}(a)_{\lambda+ y- x}\varPsi(b)]_{\lambda+ \mu+ y- x}c]+ (-1)^{|\psi||a|}[[\varPsi(a)_{\lambda}\psi_{y- x}(b)]_{\lambda+ \mu+y- x}c]\\&+ (-1)^{|\psi|(|a|+ |b|)} [[\varPsi(a)_{\lambda}b]_{\lambda+ \mu}\psi_{y- x}(c)])\\
=&([[\phi_x\psi_{y- x}(a)_{\lambda+ y}\varPhi\varPsi(b)]_{\lambda+ \mu+ y}c]\\&+ (-1)^{|\phi|(|\psi|+ |a|)}[[\varPhi\psi_{y- x}(a)_{\lambda+ y- x}\varPsi\phi_x(b)]_{\lambda+ \mu+ y}c]\\&+ (-1)^{|\phi|(|\psi|+|a|+ |b|)}[[\varPhi\psi_{y- x}(a)_{\lambda+ y- x}\varPsi (b)]_{\lambda+ \mu+ y- x}\phi_x(c)])\\
     &+(-1)^{|a||\psi|}([[\phi_x\varPsi(a)_{\lambda+ x}\varPhi\psi_{y- x}(b)]_{\lambda+ \mu+ y}c]+ (-1)^{|a||\phi|} [[\varPhi\varPsi(a)_{\lambda}\phi_x\psi_{y- x}(b)]_{\lambda+ \mu+ y}c]\\&+ (-1)^{(|a|+ |\psi|+ |b|)|\phi|}[[\varPhi\varPsi(a)_{\lambda}\psi_{y- x}(b)]_{\lambda+ \mu+ y- x}\phi_x(c)])\\
&+ (-1)^{|\psi|(|a|+ |b|)} ([[\varPsi\phi_x(a)_{\lambda+ x}\varPhi(b)]_{\lambda+ \mu+ x}\psi_{y- x}(c)]+ (-1)^{|a||\phi|}[[\varPhi\varPsi(a)_{\lambda}\phi_x(b)]_{\lambda+ \mu+ x}\psi_{y- x}(c)]\\&+ (-1)^{(|a|+ |b|)|\phi|}[[\varPhi\varPsi(a)_{\lambda}b]_{\lambda+ \mu}\phi_x\psi_{y- x}(c)])\\=&([[\phi_x\psi_{y- x}(a)_{\lambda+ y}\varPhi\varPsi(b)]_{\lambda+ \mu+ y}c]+ (-1)^{|a|(|\phi|+|\psi|)} [[\varPhi\varPsi(a)_{\lambda}\phi_x\psi_{y- x}(b)]_{\lambda+ \mu+ y}c]\\&+ (-1)^{(|a|+ |b|)(|\phi|+|\psi|)}[[\varPhi\varPsi(a)_{\lambda}b]_{\lambda+ \mu}\phi_x\psi_{y- x}(c)])\\&+ (-1)^{|\phi|(|\psi|+ |a|)}[[\varPhi\psi_{y- x}(a)_{\lambda+ y- x}\varPsi\phi_x(b)]_{\lambda+ \mu+ y}c]\\&+ (-1)^{|\phi|(|\psi|+|a|+ |b|)}[[\varPhi\psi_{y- x}(a)_{\lambda+ y- x}\varPsi (b)]_{\lambda+ \mu+ y- x}\phi_x(c)]\\&+(-1)^{|a||\psi|}[[\phi_x\varPsi(a)_{\lambda+ x}\varPhi\psi_{y- x}(b)]_{\lambda+ \mu+ y}c]\\&+(-1)^{|\phi|(|a|+ |\psi|+ |b|)+|a||\psi|}[[\varPhi\varPsi(a)_{\lambda}\psi_{y- x}(b)]_{\lambda+ \mu+ y- x}\phi_x(c)]\\&+ (-1)^{|\psi|(|a|+ |b|)} [[\varPsi\phi_x(a)_{\lambda+ x}\varPhi(b)]_{\lambda+ \mu+ x}\psi_{y- x}(c)]\\&+ (-1)^{|a||\phi|+|\psi|(|a|+ |b|)}[[\varPhi\varPsi(a)_{\lambda}\phi_x(b)]_{\lambda+ \mu+ x}\psi_{y- x}(c)].
	\end{align*}
	Similarly, we can obtain that 
	\begin{align*}
	\psi_{y- x}&\phi_x([[a_{\lambda}b]_{\lambda+ \mu}c])\\=& ([[\psi_{y- x}\phi_x(a)_{\lambda+ y}\varPsi\varPhi(b)]_{\lambda+ \mu+ y}c]
	+ (-1)^{|a|(|\psi|+ |\phi|)} [[\varPsi\varPhi(a)_{\lambda}\psi_{y- x}\phi_x(b)]_{\lambda+ \mu+ y}c]\\& +(-1)^{(|\psi|+|\phi|)(|a|+ |b|)}[[\varPsi\varPhi(a)_{\lambda}b]_{\lambda+ \mu}\psi_{y- x}\phi_x(c)])\\&+(-1)^{|\psi|(|\phi|+ |a|)}[[\varPsi\phi_x(a)_{\lambda+ x}\varPhi\psi_{y- x}(b)]_{\lambda+ \mu+ y}c]\\&+ (-1)^{|\psi|(|\phi|+|a|+ |b|)}[[\varPsi\phi_x(a)_{\lambda+x} \varPhi(b)]_{\lambda+ \mu+ x}\psi_{y- x}(c)]
	\end{align*}
	\begin{align*}
 &+(-1)^{|a||\phi|}[[\psi_{y- x}\varPhi(a)_{\lambda+y -x}\varPsi\phi_x(b)]_{\lambda+ \mu+ y}c]\\&+(-1)^{|\psi|(|a|+ |\phi|+ |b|)+|a||\phi|}[[\varPsi\varPhi(a)_{\lambda}\phi_x(b)]_{\lambda+ \mu+  x}\psi_{y- x}(c)]\\&+ (-1)^{|\phi|(|a|+ |b|)} [[\varPhi\psi_{y- x}(a)_{\lambda+y- x}\varPsi(b)]_{\lambda+ \mu+y- x}\phi_x(c)]\\&+(-1)^{|a||\psi|+|\phi|(|a|+ |b|)}[[\varPsi\varPhi(a)_{\lambda}\psi_{y- x}(b)]_{\lambda+ \mu+ y- x}\phi_x(c)]
	\end{align*}
By combining the above two equations, we have
\begin{align*}
[\phi_x\psi]_{y}([[a_{\lambda}b]_{\lambda+ \mu}c])=& [[[\phi_x\psi]_{y}(a)_{\lambda+ y}\varPsi\varPhi(b)]_{\lambda+ \mu+ y}c]+ (-1)^{|a|(|\psi|+ |\phi|)} [[\varPsi\varPhi(a)_{\lambda}[\phi_x\psi]_{y}(b)]_{\lambda+ \mu+ y}c]\\&+(-1)^{(|\psi|+|\phi|)(|a|+ |b|)}[[\varPsi\varPhi(a)_{\lambda}b]_{\lambda+ \mu}[\phi_x\psi]_{y}(c)].
\end{align*}	
This implies that $[\phi_x\psi]_{y} \in CTDer_{\varPhi\varPsi}(\mathcal{R})_{\alpha+\beta}$. 

This completes the proof.
\end{proof}

\begin{corollary}
Let $\mathcal{R}$ be a Lie conformal superalgebra and $G$ be a subgroup of $Aut(\mathcal{R})$. By assuming $\varPhi^{2}=\varPhi\in G$ such that $\varPhi$ commutes with each element of $CTDer_{\varPhi}(\mathcal{R})$. Then $CTDer_{\varPhi}(\mathcal{R})$ is a Lie conformal superalgebra.
\end{corollary}
\begin{proof} Since $CTDer_{\varPhi}(\mathcal{R})$ is a $C[\partial]$- module. So it is sufficient to prove that $CTDer_{\varPhi}(\mathcal{R})$ is closed under $\lambda$-bracket. For any $\phi_{x} \in CTDer_{\varPhi}(\mathcal{R})_\alpha$, $\psi_{x} \in CTDer_{\varPhi}(\mathcal{R})_\beta$ and homogeneous elements $a,b,c\in \mathcal{R}$, with a similar discussion in the proof of  Proposition \ref{prop5.5}, we have
	\begin{align*}
	[\phi_x\psi]_{y}([[a_{\lambda}b]_{\lambda+ \mu}c])=& [[[\phi_x\psi]_{y}(a)_{\lambda+ y}\varPhi^{2}(b)]_{\lambda+ \mu+ y}(c)]+ (-1)^{|a|(|\psi|+ |\phi|)} [[\varPhi^{2}(a)_{\lambda}[\phi_x\psi]_{y}(b)]_{\lambda+ \mu+ y}c]\\&+(-1)^{(|\psi|+|\phi|)(|a|+ |b|)}[[\varPhi^{2}(a)_{\lambda}b]_{\lambda+ \mu}[\phi_x\psi]_{y}(c)].
	\end{align*}	
	Since $\varPhi^{2}=\varPhi$, we have
	\begin{align*}
	[\phi_x\psi]_{y}([[a_{\lambda}b]_{\lambda+ \mu}c])=& [[[\phi_x\psi]_{y}(a)_{\lambda+ y}\varPhi(b)]_{\lambda+ \mu+ y}c]+ (-1)^{|a|(|\psi|+ |\phi|)} [[\varPhi(a)_{\lambda}[\phi_x\psi]_{y}(b)]_{\lambda+ \mu+ y}c]\\&+(-1)^{(|\psi|+|\phi|)(|a|+ |b|)}[[\varPhi(a)_{\lambda}b]_{\lambda+ \mu}[\phi_x\psi]_{y}(c)].
	\end{align*}
Thus, $[\phi_x\psi]_{y}\in CTDer_{\varPhi}(\mathcal{R})_{\alpha+\beta}$.
	
\end{proof}
\begin{proposition}
Let $\mathcal{R}$ be a Lie conformal superalgebra and $G$ be a subgroup of $Aut(\mathcal{R})$. Suppose that $\varPhi, \varPsi \in G$ and $\phi_{x}\in CTDer_{\varPhi}(\mathcal{R})$. Then $\varPsi\phi_{x} \in
	CTDer_{\varPsi\varPhi, \varPsi}(\mathcal{R})$ and $\phi_{x}\varPsi \in
	CTDer_{\varPhi\varPsi, \varPsi}(\mathcal{R})$
\end{proposition}
\begin{proof}
Let $\phi_{x} \in CTDer_{\varPhi}(\mathcal{R})_\alpha$, for homogeneous elements $a, b ,c \in \mathcal{R}$, we have
	\begin{align*}
	\varPsi\phi_x([[a_{\lambda}b]_{\lambda+\mu}c])&=[[\varPsi(\phi_x(a))_{\lambda+x}\varPsi\varPhi(b)]_{\lambda+\mu+x}\varPsi(c)]+(-1)^{|\phi||a|}[[\varPsi\varPhi(a)_{\lambda}\varPsi(\phi_x(b))]_{\lambda+\mu+x}\varPsi(c)]\\&+(-1)^{|\phi|(|a|+|b|)}[[\varPsi\varPhi(a)_{\lambda}\varPsi(b)]_{\lambda+\mu}\varPsi(\phi_x(c))],
	\end{align*}
Hence, $\varPsi\phi_x\in CTDer_{\varPsi\varPhi, \varPsi}(\mathcal{R})_\alpha$. Similarly, we have $\phi_x\varPsi\in CTDer_{\varPhi\varPsi, \varPsi}(\mathcal{R})_\alpha.$
	\end{proof}

\begin{proposition}
Let $\mathcal{R}$ be a Lie conformal superalgebra and $G$ be a subgroup of $Aut(\mathcal{R})$. Suppose that $\phi_{x}\in CTDer_{\varPhi}(\mathcal{R})_\alpha$ is an element such that $(\phi_x\varPhi-\varPhi\phi_x)(\mathcal{R}) \subseteq Z(\mathcal{R}).$ Then $[[{\mathcal{R}}_{\lambda}\mathcal{R}]_{\lambda+\mu} \mathcal{R}]\subseteq Ker(\phi_x\varPhi-\varPhi\phi_x)$.
	
\end{proposition}

\begin{proof}For any homogeneous elements $a, b ,c \in \mathcal{R}$ and $\phi_{x} \in CTDer_{\varPhi}(\mathcal{R})_\alpha$, we have
	\begin{align*}
	\varPhi\phi_x([[a_{\lambda}b]_{\lambda+\mu}c])&=[[\varPhi(\phi_x(a))_{\lambda+x}\varPhi^2(b)]_{\lambda+\mu+x}\varPhi(c)]+(-1)^{|\phi||a|}[[\varPhi^2(a)_{\lambda}\varPhi(\phi_x(b))]_{\lambda+\mu+x}\varPhi(c)]\\&+(-1)^{|\phi|(|a|+|b|)}[[\varPhi^2(a)_{\lambda}\varPhi(b)]_{\lambda+\mu}\varPhi(\phi_x(c))],
	\end{align*}
	and
		\begin{align*}
	\phi_x\varPhi([[a_{\lambda}b]_{\lambda+\mu}c])&=[[\phi_x(\varPhi(a))_{\lambda+x}\varPhi^2(b)]_{\lambda+\mu+x}\varPhi(c)]+(-1)^{|\phi||a|}[[\varPhi^2(a)_{\lambda}\phi_x(\varPhi(b))]_{\lambda+\mu+x}\varPhi(c)]\\&+(-1)^{|\phi|(|a|+|b|)}[[\varPhi^2(a)_{\lambda}\varPhi(b)]_{\lambda+\mu}\phi_x(\varPhi(c))],
	\end{align*}
	Since $(\phi_x\varPhi-\varPhi\phi_x)(\mathcal{R}) \subseteq Z(\mathcal{R})$, we have
	\begin{align}
	(\phi_x\varPhi-\varPhi\phi_x)([[a_{\lambda}b]_{\lambda+\mu}c])= 0
	\end{align}
	Consequently, $[[{\mathcal{R}}_{\lambda}\mathcal{R}]_{\lambda+\mu} \mathcal{R}]\subseteq Ker(\phi_x\varPhi-\varPhi\phi_x)$.
\end{proof}
\par Later in this section we will assume that  $\mathcal{R}$ is a Lie conformal superalgebra. We define and classify all possible \emph{$(A,B,C,D)$-derivation} of degree $\alpha$ of $\mathcal{R}$.
\begin{definition}\label{def6.1}
	\begin{em}
		Let $\mathcal{R}$ be a Lie conformal superalgebra. A conformal linear map $\phi_{x} \in gc(\mathcal{R})$ is called an \emph{(A,B,C,D)-derivation} of degree $\alpha$ of $\mathcal{R}$ if there exists $A,B,C,D\in\mathbb{C}$ such that
		\begin{align}\label{eq12}
		A\phi_x([[a_{\lambda}b]_{\lambda+\mu}c])=&B[[\phi_x(a)_{\lambda+x}b]_{\lambda+\mu+x}c]+(-1)^{|\phi||a|}C[[a_{\lambda}\phi_x(b)]_{\lambda+\mu+x}c]\nonumber\\&+(-1)^{|\phi|(|a|+|b|)}D[[a_{\lambda}b]_{\lambda+\mu}\phi_x(c)],
		\end{align}
		for any homogeneous elements $a,b,c \in \mathcal{R}$.
	\end{em}
\end{definition}
The set of all $(A,B,C,D)$-derivations of $\mathcal{R}$ is denoted by $\varPhi(A,B,C,D)$. Obviously, every conformal triple derivation on  $\mathcal{R}$ is in fact $(1,1,1,1)$-derivation of Lie conformal superalgebra $\mathcal{R}$.

\begin{proposition}\label{prop6.2}
	For all $A,B,C,D \in \mathbb{C}$, the set $\varPhi(A,B,C,D)$ can be written as
	\begin{equation}\label{eq6.2}
	\varPhi(A, B, C, D)= \varPhi(2A, B+ C, B+ C, 2D)\cap \varPhi(0, B- C, C- B, 0).
	\end{equation}
	Moreover, we have $\varPhi(A, B, C, D)= \varPhi(A, C, B, D) =\varPhi(tA, tB, tC, tD)$, for any $t \in \mathbb{C} / \{0\}$.
\end{proposition}
\begin{proof}
	For any $A,B,C,D \in \mathbb{C}$ and $\phi_{x} \in \varPhi(A,B,C,D)$, we have
	\begin{eqnarray}\nonumber A\phi_x([[a_{\lambda}b]_{\lambda+\mu}c])&=B[[\phi_x(a)_{\lambda+x}b]_{\lambda+\mu+x}c]+(-1)^{|\phi||a|}C[[a_{\lambda}\phi_x(b)]_{\lambda+\mu+x}c]\\&\label{eq13}+(-1)^{|\phi|(|a|+|b|)}D[[a_{\lambda}b]_{\lambda+\mu}\phi_x(c)],\end{eqnarray} and 
	\begin{eqnarray}\nonumber A\phi_x([[b_{\lambda}a]_{\lambda+\mu}c])&=B[[\phi_x(b)_{\lambda+x}a]_{\lambda+\mu+x}c]+(-1)^{|\phi||b|}C[[b_{\lambda}\phi_x(a)]_{\lambda+\mu+x}c]\\&\label{eq14}+(-1)^{|\phi|(|b|+|a|)}D[[b_{\lambda}a]_{\lambda+\mu}\phi_x(c)],
	\end{eqnarray}
	where $a,b,c \in \mathcal{R}$ are homogeneous elements .
	By the conformal skew-symmetry property of $\mathcal{R}$, (\ref{eq14}) can be transform as
	\begin{align*}
	-(-1)^{|a||b|}A\phi_x([[a_{-\partial-\lambda}b]_{\lambda+\mu}c])=&-(-1)^{(|\phi|+|b|)|a|}B[[a_{-\partial-\lambda-x}\phi_x(b)]_{\lambda+\mu+x}c]
	\\-&(-1)^{|\phi||b|+(|\phi|+|a|)|b|}C[[\phi_x(a)_{-\partial-\lambda} b]_{\lambda+\mu+x}c]
	\\-&(-1)^{|\phi|(|b|+|a|)+|a||b|}D[[a_{-\partial-\lambda}b]_{\lambda+\mu}\phi_x(c)].
	\end{align*}
	Hence, we have
	\begin{eqnarray}
	A\phi_x([[a_{-\partial-\lambda}b]_{\lambda+\mu}c])=&C[[\phi_x(a)_{-\partial-\lambda} b]_{\lambda+\mu+x}c]+(-1)^{|\phi||a|}B[[a_{-\partial-\lambda-x}\phi_x(b)]_{\lambda+\mu+x}c]\nonumber\\&\label{eq16}+(-1)^{|\phi|(|b|+|a|)}D[[a_{-\partial-\lambda}b]_{\lambda+\mu}\phi_x(c)].\label{eq6.5}
	\end{eqnarray}
	By the conformal sesquilinearity and swapping the position of $\lambda$ and $\mu$ in (\ref{eq6.5}), we can obtain that
	
	\begin{eqnarray}\nonumber A\phi_x([[a_{\lambda}b]_{\lambda+\mu}c])&=C[[\phi_x(a)_{\lambda+x}b]_{\lambda+\mu+x}c]+(-1)^{|\phi||a|}B[[a_{\lambda}\phi_x(b)]_{\lambda+\mu+x}c]\\&\label{eq6.6}+(-1)^{|\phi|(|a|+|b|)}D[[a_{\lambda}b]_{\lambda+\mu}\phi_x(c)].\end{eqnarray}
	By (\ref{eq13}) and (\ref{eq6.6}) we have
	\begin{eqnarray}\nonumber
	2A\phi_x([[a_{\lambda}b]_{\lambda+\mu}c])&=(B+C)[[\phi_x(a)_{\lambda+x}b]_{\lambda+\mu+x}c]+(-1)^{|\phi||a|}(C+B)[[a_{\lambda}\phi_x(b)]_{\lambda+\mu+x}c]\\&\label{eq17}+(-1)^{|\phi|(|a|+|b|)}2D[[a_{\lambda}b]_{\lambda+\mu}\phi_x(c)],\end{eqnarray}and
	\begin{eqnarray}\label{eq18}
	0=(B-C)[[\phi_x(a)_{\lambda+x}b]_{\lambda+\mu+x}c]+(-1)^{|\phi||a|}(C-B)[[a_{\lambda}\phi_x(b)]_{\lambda+\mu+x}c]+0.\end{eqnarray}
	Thus,
	\nonumber\begin{align}\phi_{x}\in \varPhi(2A, B+ C, B+ C, 2D)\cap \varPhi(0, B- C, C- B, 0).\end{align}
	This implies that \begin{eqnarray}\label{eq20}\varPhi(A,B,C,D)\subseteq \varPhi(2A, B+ C, B+ C, 2D)\cap \varPhi(0, B- C, C- B, 0)\end{eqnarray}
	
	Conversely, it is not difficult to check that 
	\begin{eqnarray}\label{eq202}\varPhi(2A, B+ C, B+ C, 2D)\cap \varPhi(0, B- C, C- B, 0)\subseteq \varPhi(A,B,C,D)\end{eqnarray}
	by using (\ref{eq17}) and (\ref{eq18}).
	Thus, we have 
	\begin{align}\label{eq22}
	\varPhi(A, B, C, D)= \varPhi(2A, B+ C, B+ C, 2D)\cap \varPhi(0, B- C, C- B, 0). \end{align}	
	This completes the proof.\end{proof}

\begin{proposition}\label{prop6.4}
	For any $A,B,C,D\in \mathbb{C}$ and $B\ne \pm C$,  we have 
	\begin{align}
	\varPhi(A, B, C, D)= \varPhi(\frac{A}{B+ C}, 0, 1, \frac{D}{B+ C})=\varPhi(\frac{A}{B+ C}, 1, 0, \frac{D}{B+ C}).
	\end{align}	
\end{proposition}
\begin{proof}By equation (\ref{eq6.2}), we have \nonumber\begin{align}\varPhi(A, B, C, D)&= \varPhi(2A, B+ C, B+ C, 2D)\cap \varPhi(0,B- C, C- B, 0)\\&= \varPhi(\frac{2A}{B+ C}, 1, 1, \frac{2D}{B+C})\cap \varPhi(0, -1, 1, 0)\\&= \varPhi(\frac{2A}{B+ C}, 0, 2, \frac{2D}{B+ C})\\&= \varPhi(\frac{A}{B+ C} ,0 ,1 , \frac{D}{B+ C}).	\end{align}
	This completes the proof.
\end{proof}

By combining the above mentioned fact with the results of Proposition \ref{prop6.2} and \ref{prop6.4}, we get the  following set of identities:

\begin{enumerate}
	\item For $B+C=0$, either $B=C=0$ or $B=-C\neq0$, then
\begin{eqnarray}\label{eq28}\varPhi(A,B,C,D)=\left\{\begin{array}{l} \varPhi(A,0,0,D),\qquad D\neq 0, B=C=0,\\
\varPhi(A,0,0,0), \qquad D=0, B=C=0,\\
\varPhi(A,B,-B,D),\qquad D\neq 0 , B=-C\neq 0,\\
\varPhi(A,B,-B,0), \qquad D=0,  B=-C\neq 0.\end{array}\right.\end{eqnarray}	
	
\item For $B+C\neq0$, either $B-C\neq0$ or $B=C\neq 0$,	then
	\begin{eqnarray}\label{eq29}\varPhi(A,B,C,D)=\left\{\begin{array}{l} \varPhi(A,B,C,D),\qquad D\neq 0, B-C\neq0,\\
	\varPhi(A,B,C,0), \qquad D=0, B-C\neq0,\\
	\varPhi(A,B,B,D),\qquad D\neq 0 , B=C\neq 0,\\
	\varPhi(A,B,B,0), \qquad D=0,  B=C\neq 0.\end{array}\right.\end{eqnarray}
\end{enumerate}

According to Definition \ref{def6.1}, if $\frac{B}{A} =\frac{C}{A} =\frac{D}{A} = \delta$, then  $(A, B, C, D)$-derivation is a $\delta$-derivation (a kind of quasi- triple derivation) of $\mathcal{R}$. Now equations (\ref{eq28}) and (\ref{eq29}) changed to,

	\begin{eqnarray}\label{eq30}
	\varPhi(A,B,C,D)=\left\{\begin{array}{l} \varPhi(\delta,0,0,1),\qquad D\neq 0, B=C=0,\\
\varPhi(\delta,0,0,0), \qquad D=0, B=C=0,\\
\varPhi(\delta,1,-1,1),\qquad D\neq 0 , B=-C\neq 0,\\
\varPhi(\delta ,1, -1, 0), \qquad D=0,  B=-C\neq 0,\end{array}\right.\end{eqnarray}
and
	\begin{eqnarray}\label{eq31}\varPhi(A,B,C,D)=\left\{\begin{array}{l} \varPhi(\delta,1,0,\delta^{'}),\qquad D\neq 0, B-C\neq0,\\
\varPhi(\delta,1,0,0), \qquad D=0, B-C\neq0,\\
\varPhi(\delta,1,1,\delta^{'}),\qquad D\neq 0 , B=C\neq 0,\\
\varPhi(\delta ,1, 1, 0), \qquad D=0,  B=C\neq 0.\end{array}\right.\end{eqnarray}

In order to describe the above situation more clearly, we get various results as listed below coincides with the results in \cite{Zhao-Ma-Chen}.

\begin{enumerate}
	\item $\varPhi(\delta,0,0,1)= \{\phi_{x} \in gc(\mathcal{R})|\delta \phi_{x}([[a_\lambda b]_{\lambda+\mu} c]) = (-1)^{|\phi|(|a|+|b|)}[[a_\lambda b]_{\lambda+\mu}\phi_{x}(c)],\\ \forall a,b,c \in \mathcal{R} \},$\\
	\item $\varPhi(\delta,0,0,0)= \{\phi_{x}\in gc(\mathcal{R})|\delta([[a_{\lambda}b]_{\lambda+ \mu}c]) =0, \forall a,b,c\in \mathcal{R} \},$\\
	\item $\varPhi(\delta,1,-1,1)= \varPhi(\delta,0,0,1)\cap\varPhi(0,-1,1,0),$\\
	\item $\varPhi(\delta, 1, -1, 0)= \varPhi(\delta, 0, 0, 0)\cap \varPhi(0, -1, 1, 0),$\\
	\item $\varPhi(\delta,1,0,\delta^{'})=\varPhi(2\delta,1,1,2\delta^{'})\cap \varPhi(0,1,-1,0),$\\
	\item $\varPhi(\delta,1,0,0)=  \varPhi(2\delta , 1, 1, 0)\cap \varPhi(0, 1, -1, 0),$\\
   \item $\varPhi(\delta,1,1,\delta^{'})= \{\phi_{x}\in gc(\mathcal{R})| \delta \phi_{x}([[a_\lambda b]_{\lambda+\mu} c])= [[\phi_{x}(a)_{\lambda+x}b]_{\lambda+x + \mu}c]\\+(-1)^{|\phi||a|}[[a_{\lambda} \phi_{x}(b)]_{\lambda+x + \mu}c]+ (-1)^{|\phi|(|a|+|b|)}\delta^{'} [[a_{\lambda}b]_{\lambda+ \mu}\phi_{x}(c)], \forall a,b,c\in \mathcal{R} \},$\\
\item $\varPhi(\delta ,1, 1, 0)= \{\phi_{x} \in gc(\mathcal{R})| \delta \phi_{x}([[a_{\lambda}b]_{\lambda+ \mu} c] )= [[\phi_{x}(a)_{\lambda+x} b]_{\lambda+x + \mu} c]\\ + (-1)^{|\phi||a|}[[a_{\lambda}\phi_{x}(b)]_{\lambda+x + \mu}c], \forall a,b,c \in \mathcal{ R}\}.$\end{enumerate}

\begin{remark}
In particular, we can obtain some connections between $(A,B,C,D)$-derivations and conformal triple derivations on Lie conformal superalgebra as follows.
\begin{enumerate}
       \item $\varPhi(0,0,0,0)= gc(\mathcal{R}),$\\
	\item $\varPhi(1,1,0,0)=TC(\mathcal{R}),$\\
	\item $\varPhi(0,1,0,-1)=TQC(\mathcal{R}),$\\
	\item $\varPhi(1,1,1,1)= CTDer(\mathcal{R}).$
\end{enumerate}
\end{remark}

\begin{theorem}
Let $\mathcal{R}$ be a Lie conformal superalgebra. The $\varPhi(A,B,C,D)$-derivation of degree $\alpha$ of $\mathcal{R}$ is isomorphic to one of the following cases
\begin{align*}
&(1)\varPhi(\delta,0,0,1),\quad (2)\varPhi(\delta,0,0,0),\quad (3)\varPhi(\delta,1,-1,1),\quad(4)\varPhi(\delta, 1, -1, 0),\\
&(5)\varPhi(\delta,1,0,\delta^{'}),\quad(6)\varPhi(\delta,1,0,0),\quad(7)\varPhi(\delta,1,1,\delta^{'}),\quad(8)\varPhi(\delta ,1, 1, 0),
\end{align*}
where $\delta, \delta^{'} \in \mathbb{C}.$
\end{theorem}

\section{Classification of conformal triple homomorphisms on simple Lie conformal superalgebras}

In this section, we first introduce the definition of triple homomorphism of a Lie conformal superalgebra. Then, the conformal triple homomorphism of all finite simple Lie conformal superalgebras is also characterized.

\begin{definition}\label{def5.1}
	\begin{em}
		Let $\mathcal{A}, \mathcal{B}$ be two Lie conformal superalgebras. A $\mathbb{C}[\partial]$-module homomorphism $f: \mathcal{A} \to \mathcal{B}$ is called:
		\begin{enumerate}
			\item a \emph{homomorphism} if it satisfies $f([a_\lambda b])=[f(a)_\lambda f(b)]$ for any $a, b \in \mathcal{A}$.
			\item a \emph{anti-homomorphism} if it satisfies $f([a_\lambda b])=-[f(a)_\lambda f(b)]$ for any $a,b \in \mathcal{A}$.
			\item a \emph{triple homomorphism} if it satisfies $f([a_\lambda [b_\mu c]])=[f(a)_\lambda [f(b)_\mu f(c)]]$ for any $a,b,c \in \mathcal{A}$.
		\end{enumerate}
	\end{em}		
\end{definition}
\begin{remark}\label{rm5.2}
	\begin{em}
		Obviously, due to the skew-symmetry of Lie conformal superalgebra, the following two equations are equivalent:
		\begin{enumerate}
			\item $f([a_\lambda [b_\mu c]])=[f(a)_\lambda [f(b)_\mu f(c)]]$.
			\item $f([[a_\lambda b]_{\lambda+\mu} c])=[[f(a)_\lambda f(b)]_{\lambda+\mu} f(c)]$.
		\end{enumerate}
	\end{em}
\end{remark}

Let $\mathcal{A}$ be a Lie conformal superalgebra. For a subset $S$ of $\mathcal{A}$, the \emph{enveloping Lie conformal superalgebra} of $S$ is by definition the Lie conformal subalgebra of $\mathcal{A}$ generated by $S$. A Lie conformal superalgebra is called \emph{indecomposable} if it cannot be written as a direct sum of two nontrivial ideals.

\begin{definition}\label{def5.2}
	\begin{em}
		Let $\mathcal{A}, \mathcal{B}$ be two Lie conformal superalgebras. A $\mathbb{C}[\partial]$-module homomorphism $f: \mathcal{A} \to \mathcal{B}$ is called a \emph{direct sum} of $f_1$ and $f_2$, if $f=f_1+f_2$ and there exists ideals $I_1$, $I_2$ of the enveloping Lie conformal superalgebra of $f(\mathcal{A})$ such that $I_1\cap I_2=0$ and $f_1(\mathcal{A}) \subseteq I_1$, $f_2(\mathcal{A}) \subseteq I_2$. 
	\end{em}		
\end{definition}
\begin{remark}
	\begin{em}	
		Obviously, both homomorphisms and anti-homomorphisms of Lie conformal superalgebra are triple homomorphisms, but the converse does not always hold. Moreover, the sum of a homomorphism and an anti-homomorphism is also a triple homomorphism. Actually, there exists some other connection between them which is not just inclusion relationship. 
	\end{em}	
\end{remark}
In the sequel, we always use $E$ to denote the enveloping Lie conformal superalgebra of $f(\mathcal{A})$. Moreover, we suppose that $\mathcal{A}$ is a finite simple Lie conformal superalgebra, and $E$ is centerless and can be decomposed into a direct sum of indecomposable ideals.
\begin{lemma}\label{lm5.4}
	Suppose that $\mathcal{A}, \mathcal{B}$ are Lie conformal superalgebras and $f$ is a triple homomorphism from $\mathcal{A}$ to $\mathcal{B}$. Then there exists a conformal homomorphism $\delta_f: \mathcal{A} \to \mathcal{B}$ such that for all $a \in \mathcal{A}$ one has $fada_x=ad \delta_f(a)_xf$.
\end{lemma}
\begin{proof}
	Similar to the discussion of the proof Lemma \ref{lm4.1}, there exists some finite index set $I, J \subseteq \mathbb{Z^+}$ and $a= \sum_{i\in I, j\in J}\frac{d^i}{d\lambda^i}[{a_{ij}^1}_\lambda {a_{ij}^2}]|_{\lambda=0}$. Define a $\mathbb{C}$-linear map $\delta_f: \mathcal{A} \to \mathcal{B}$, such that for any $a\in \mathcal{A}$, $\delta_f(a)= \sum_{i\in I, j\in J}\frac{d^i}{d\lambda^i}[f(a_{ij}^1)_{\lambda} f({a_{ij}^2})]|_{\lambda=0}$.

	At first, it is sufficient to prove that  $\delta_f(a)$ is independent of the expression of $a$. For proving it, take 
	\begin{align*}
	A= \sum_{i\in I, j\in J}\frac{d^i}{d\lambda^i}[f(a_{ij}^1)_{\lambda} f({a_{ij}^2})]|_{\lambda=0}, \quad B=\sum_{s\in I^{'}, t\in J^{'}}\frac{d^s}{d\lambda^s}[f(b_{st}^1)_{\lambda} f({b_{st}^2})]|_{\lambda=0},
	\end{align*}
	where $a$ can also be expressed in the form $a= \sum_{s\in I^{'}, t\in J^{'}}\frac{d^s}{d\lambda^s}[{b_{st}^1}_{\lambda} {b_{st}^2}]|_{\lambda=0}$. Since $f$ is a conformal triple homomorphism from $\mathcal{A}$ to $\mathcal{B}$, for any $c \in \mathcal{A}$, we can deduce that 
	\begin{align*}
	[f(c)_\mu (A- B)]&= [f(c)_\mu \sum_{i\in I, j\in J}\frac{d^i}{d\lambda^i}[f(a_{ij}^1)_{\lambda} f({a_{ij}^2})]|_{\lambda=0}]]- [f(c)_\mu \sum_{s\in I^{'}, t\in J^{'}}\frac{d^s}{d\lambda^s}[f(b_{st}^1)_{\lambda} f({b_{st}^2})]|_{\lambda=0}]]\\
	&= \sum_{i\in I, j\in J}\frac{d^i}{d\lambda^i}([f(c)_\mu [f(a_{ij}^1)_{\lambda} f({a_{ij}^2})]])|_{\lambda=0}- \sum_{s\in I^{'}, t\in J^{'}}\frac{d^s}{d\lambda^s}([f(c)_\mu [f(b_{st}^1)_{\lambda} f({b_{st}^2})]])|_{\lambda=0}\\
	&= \sum_{i\in I, j\in J}\frac{d^i}{d\lambda^i}(f([c_\mu [{a_{ij}^1}_{\lambda} {a_{ij}^2}]]))|_{\lambda=0}- \sum_{s\in I^{'}, t\in J^{'}}\frac{d^s}{d\lambda^s}(f([c_\mu [{b_{st}^1}_{\lambda} {b_{st}^2}]]))|_{\lambda=0}\\
	&= f([c_\mu \sum_{i\in I, j\in J}\frac{d^i}{d\lambda^i}([{a_{ij}^1}_{\lambda} {a_{ij}^2}])|_{\lambda=0}])- f([c_\mu \sum_{s\in I^{'}, t\in J^{'}}\frac{d^s}{d\lambda^s}([{b_{st}^1}_{\lambda} {b_{st}^2}])|_{\lambda=0}])\\
	&= f([c_\mu a])- f([c_\mu a])\\
	&= 0.
	\end{align*}
	
	Hence, $A-B \in Z(E)$. By the assumption, we have $Z(E)=0$. Thus, we can obtain that $A=B$. Therefore, $\delta_f$ is well-defined. 
	
	Next, we need to prove that $fada_x=ad \delta_f(a)_xf$ for each $a\in \mathcal{A}$. For any $a,b \in \mathcal{A}$, where $a= \sum_{i\in I, j\in J}\frac{d^i}{d\lambda^i}[{a_{ij}^1}_\lambda {a_{ij}^2}]|_{\lambda=0}$, we have 
		\begin{align*}
	fada_x(b)&= f([a_x b])\\
	&= f([(\sum_{i\in I, j\in J}\frac{d^i}{d\lambda^i}[{a_{ij}^1}_\lambda {a_{ij}^2}]|_{\lambda=0})_x b])\\
	&= \sum_{i\in I, j\in J}\frac{d^i}{d\lambda^i}[[f({a_{ij}^1})_\lambda f({a_{ij}^2})]_x f(b)]|_{\lambda=0}\\
	&= [(\sum_{i\in I, j\in J}\frac{d^i}{d\lambda^i}[f({a_{ij}^1})_\lambda f({a_{ij}^2})]|_{\lambda=0})_x f(b)]\\
	&= [\delta_f(a)_x f(b)]\\
	&= ad \delta_f(a)_x f(b).
	\end{align*}
	Therefore, by arbitrariness of $b\in \mathcal{A}$, for each $a\in \mathcal{A}$, we can deduce that $fada_x=ad \delta_f(a)_xf$.
	
	Finally, we need to prove that $\delta_{f}$ is a conformal homomorphism, i.e. $\delta_{f}(\partial a)=\partial\delta_{f}(a)$ and $\delta_f([a_\lambda b])=[\delta_f(a)_\lambda \delta_f(b)]$ for any $a,b \in \mathcal{A}$. For any $a= \sum_{i\in I, j\in J}\frac{d^i}{d\lambda^i}[{a_{ij}^1}_\lambda {a_{ij}^2}]|_{\lambda=0} \in \mathcal{A}$, we have $\partial a= \sum_{i\in I, j\in J}\frac{d^i}{d\lambda^i}([{\partial a_{ij}^1}_\lambda {a_{ij}^2}] + [{a_{ij}^1}_\lambda \partial {a_{ij}^2}])|_{\lambda=0}$. Thus, we can obtain that
	\begin{align*}
	\delta_{f}(\partial a)&= \sum_{i\in I, j\in J}\frac{d^i}{d\lambda^i}([f(\partial a_{ij}^1)_{\lambda} f({a_{ij}^2})]+ [f({a_{ij}^1})_\lambda f({\partial a_{ij}^2})])|_{\lambda=0}\\
	&= \sum_{i\in I, j\in J}\frac{d^i}{d\lambda^i}(-\lambda [f(a_{ij}^1)_{\lambda} f({a_{ij}^2})]+ (\partial + \lambda)[f({a_{ij}^1})_\lambda f(a_{ij}^2)])|_{\lambda=0}\\
	&= \partial\sum_{i\in I, j\in J}\frac{d^i}{d\lambda^i}([f(a_{ij}^1)_{\lambda} f({a_{ij}^2})])|_{\lambda=0}\\
	&= \partial\delta_{f}(a).
	\end{align*}
	Thus, $\delta_{f}$ is a $\mathbb{C}[\partial]$-module homomorphism. For any $a,b,c \in \mathcal{A}$, we can obtain that
	\begin{align*}
	[&(\delta_{f}([a_\lambda b])-[\delta_{f}(a)_\lambda \delta_f(b)])_{\lambda+\mu} f(c)]\\
	&=[\delta_{f}([a_\lambda b])_{\lambda+\mu} f(c)]-[\delta_{f}(a)_{\lambda} [\delta_f(b)_{\mu} f(c)]]+(-1)^{|a||b|}[\delta_f(b)_{\mu} [\delta_{f}(a)_\lambda f(c)]]\\
	&=[[f(a)_\lambda f(b)]_{\lambda+\mu} f(c)]-[\delta_{f}(a)_{\lambda} (ad \delta_f(b)_\mu f(c))]+(-1)^{|a||b|}[\delta_f(b)_{\mu}(ad \delta_f(a)_\lambda f(c))]\\
	&=[[f(a)_\lambda f(b)]_{\lambda+\mu} f(c)]-[\delta_{f}(a)_{\lambda} (fadb_\mu c)]+(-1)^{|a||b|}[\delta_f(b)_{\mu}(fada_\lambda c)]\\
	&=f([[a_\lambda b]_{\lambda+\mu} c])-[\delta_{f}(a)_{\lambda} f([b_\mu c])]+(-1)^{|a||b|}[\delta_f(b)_{\mu}f([a_\lambda c])]\\
	&=f([[a_\lambda b]_{\lambda+\mu} c])-ad\delta_{f}(a)_{\lambda} f([b_\mu c])+(-1)^{|a||b|}ad\delta_f(b)_{\mu}f([a_\lambda c])\\
	&=f([[a_\lambda b]_{\lambda+\mu} c])-fada_{\lambda}([b_\mu c])+(-1)^{|a||b|}fadb_{\mu}([a_\lambda c])\\
	&=f([[a_\lambda b]_{\lambda+\mu} c])-f([a_{\lambda}[b_\mu c]])+(-1)^{|a||b|}f([b_{\mu}[a_\lambda c]])\\
	&=f([[a_\lambda b]_{\lambda+\mu} c]-[a_{\lambda}[b_\mu c]]+(-1)^{|a||b|}[b_{\mu}[a_\lambda c]])\\
	&=0.
	\end{align*}
	Due to the arbitrariness of $c$, we have $\delta_{f}([a_\lambda b])-[\delta_{f}(a)_\lambda \delta_f(b)] \in Z(E)$. Thus, $\delta_{f}([a_\lambda b])=[\delta_{f}(a)_\lambda \delta_f(b)]$ since $E$ is centerless. Hence, $\delta_{f}$ is a conformal homomorphism by the arbitrariness of $a,b\in \mathcal{A}$.
	
	This completes the proof.
\end{proof}

There are some good properties on $E$ as shown in the following lemma.
\begin{lemma}\label{lm5.6}
	Denote $E^+=Im(f+\delta_{f}),  E^-=Im(f-\delta_{f})$. Then we can obtain some results as follows: 
	\begin{enumerate}
		\item $E^+$ and $E^-$ are both ideals of $E$. 
		\item $[{E^+}_\lambda E^-]=0$.
		\item $E^+\cap E^-=0$. 
	\end{enumerate}
	In fact, $E$ can be decomposed into a direct sum of ideals $E^+$ and $E^-$.
\end{lemma}
\begin{proof}
	(1) Obviously, we have $E^+$, $E^-\subseteq E$. For any $a,b \in \mathcal{A}$, we have
	\begin{align*}
	[(f(a)+\delta_f(a))_\lambda f(b)]&=[f(a)_\lambda f(b)]+ ad \delta_f(a)_\lambda f(b)\\
	&=\delta_f([a_\lambda b])+fada_\lambda (b)\\
	&=\delta_f([a_\lambda b])+f([a_\lambda b])\\
	&=(f+\delta_{f})([a_\lambda b]).
	\end{align*}
	Thus, $E^+$ is an ideal of $E$. Similarly, $E^-$ is also an ideal of $E$.
	
	(2) For any $a,b,c \in \mathcal{A}$, we can obtain that
	\begin{align*}
	[[(f(a)&+\delta_f(a))_\lambda (f(b)-\delta_f(b))]_{\lambda+\mu}f(c)]\\
	=&[[f(a)_\lambda f(b)]_{\lambda+\mu}f(c)]-[[f(a)_\lambda \delta_f(b)]_{\lambda+\mu}f(c)]\\
	&+[[\delta_f(a)_\lambda f(b)]_{\lambda+\mu}f(c)]-[[\delta_f(a)_\lambda \delta_f(b)]_{\lambda+\mu}f(c)]\\
	=&f([[a_\lambda b]_{\lambda+\mu}c])+(-1)^{|a||b|}[[\delta_f(b)_{-\partial-\lambda} f(a)]_{\lambda+\mu}f(c)]+[[\delta_f(a)_\lambda f(b)]_{\lambda+\mu}f(c)]\\
	&-[\delta_f(a)_\lambda [\delta_f(b)_{\mu}f(c)]]+(-1)^{|a||b|}[\delta_f(b)_{\mu}[\delta_f(a)_\lambda f(c)]]\\
	=&f([[a_\lambda b]_{\lambda+\mu}c])+(-1)^{|a||b|}[(ad\delta_f(b)_{-\partial-\lambda} f(a))_{\lambda+\mu}f(c)]+[(ad \delta_f(a)_\lambda f(b))_{\lambda+\mu}f(c)]\\
	&-[\delta_f(a)_\lambda (ad\delta_f(b)_{\mu}f(c))]+(-1)^{|a||b|}[\delta_f(b)_{\mu}(ad\delta_f(a)_\lambda f(c))]\\
	=&f([[a_\lambda b]_{\lambda+\mu}c])+(-1)^{|a||b|}[(fadb_{-\partial-\lambda}(a))_{\lambda+\mu}f(c)]+[(fada_\lambda b)_{\lambda+\mu}f(c)]\\
	&-[\delta_f(a)_\lambda (fadb_{\mu}c)]+(-1)^{|a||b|}[\delta_f(b)_{\mu}(fada_\lambda c)]\\
	=&f([[a_\lambda b]_{\lambda+\mu}c])+(-1)^{|a||b|}[f([b_{-\partial-\lambda}a])_{\lambda+\mu}f(c)]+[f([a_\lambda b])_{\lambda+\mu}f(c)]\\
	&-ad\delta_f(a)_\lambda f([b_{\mu}c])+(-1)^{|a||b|}ad \delta_f(b)_{\mu}f([a_\lambda c])\\=&f([[a_\lambda b]_{\lambda+\mu}c])-[f([a_\lambda b])_{\lambda+\mu}f(c)]+[f([a_\lambda b])_{\lambda+\mu}f(c)]\\
	&-fada_\lambda ([b_{\mu}c])+(-1)^{|a||b|}fadb_{\mu}([a_\lambda c])\\
	=&f([[a_\lambda b]_{\lambda+\mu}c])-f([a_\lambda [b_{\mu}c]])+(-1)^{|a||b|}f([b_{\mu}[a_\lambda c]])\\
	=&f([[a_\lambda b]_{\lambda+\mu}c]-[a_\lambda [b_{\mu}c]]+(-1)^{|a||b|}[b_{\mu}[a_\lambda c]])\\
	=&0.
	\end{align*}
	Therefore, $[(f(a)+\delta_f(a))_\lambda (f(b)-\delta_f(b))] \in Z(E)$. Since $Z(E)=0$, we obtain that $[(f(a)+\delta_f(a))_\lambda (f(b)-\delta_f(b))]=0$. 
	
	(3) For any $a \in E^+\cap E^-$, we have $[a_\lambda E^+]=[a_\lambda E^-]=0$ by (2). Thus, for any $b \in \mathcal{A}$, we have $[a_\lambda (f(b)+\delta_f(b))]=[a_\lambda (f(b)-\delta_f(b))]=0$. Then, $[a_\lambda f(b)]=0$, i.e., $a \in Z(E)$. Since $Z(E)=0$, we can deduce that $a=0$.
	
	This completes the proof.
\end{proof}

We first consider the simplest situation of $E$ to investigate the conformal triple homomorphism from $\mathcal{A}$ to $\mathcal{B}$. 
\begin{lemma}\label{lm5.7}
	If $E$ is indecomposable, then $f$ is either a conformal homomorphism or a conformal anti-homomorphism from $\mathcal{A}$ to $\mathcal{B}$.
\end{lemma}
\begin{proof}
	For any $a \in \mathcal{A}$, set $a^+=\frac{1}{2}(f(a)+\delta_f(a)), \ a^-=\frac{1}{2}(f(a)-\delta_f(a))$. Obviously, we have $a^+ \in E^+, a^-\in E^-$ and $f(a)=a^++a^-$. Therefore, $E \subseteq E^+ +  E^-$. By Lemma \ref{lm5.6}, we have $E = E^+ \oplus  E^-$. Since $E$ is indecomposable, then either $E^+$ or $E^-$ must be trivial. If $E^+$ is trivial, i.e., $(f+\delta_f)([a_{\lambda}b])=0$, then $f([a_\lambda b])=-\delta_f([a_\lambda b])=-[f(a)_{\lambda}f(b)].$ Thus, $f$ is an anti-homomorphism. Similarly, if $E^-$ is trivial, then $f$ is a homomorphism. 
	
	This completes the proof.
\end{proof}

Finally, we put forward the main result in this section by the following theorem.
\begin{theorem}
	Assume that $\mathcal{A}$ is a finite simple Lie conformal superalgebra, and $E$ is centerless and can be decomposed into a direct sum of indecomposable ideals. Then $f$ is either a homomorphism or an anti-homomorphism or a direct sum of a homomorphism and an anti-homomorphism from $\mathcal{A}$ to $\mathcal{B}$.
\end{theorem}
\begin{proof}
	Obviously, it remains to prove the Theorem \ref{lm5.7} in case $E$ is decomposable. By the assumption, $E$ can be written as the sum $E=E_1\oplus E_2\oplus ...\oplus E_n$, where each $E_i$ is an indecomposable ideal of $E$. Since $E$ is centerless, each $E_i$ is also centerless by Lemma \ref{lm5.6}, 
	\par Let $p_i$ be the projection of $E$ into $E_i$. Then we have $f=\sum_{i=1}^n p_if$ and $p_if$ is a triple homomorphism from $\mathcal{A}$ to $E_i$, and $E_i$ is the enveloping Lie conformal superalgebras of $p_if(\mathcal{A})$  for $i= 1,2,...,n$. Since each $E_i$ is indecomposable, by Lemma \ref{lm5.7}, we can deduce that $p_if$ is either a homomorphism or an anti-homomorphism from $\mathcal{A}$ to $E_i$. Set $I=\{i~~|~~p_if \ is\  a\  homomorphism.\}$ and $J=\{j~~|~~ p_j f \ is\  an\  anti-homomorphism.\}$. Set $E_I=\sum_{i \in I} \oplus E_i$ and $E_J=\sum_{j \in J}\oplus E_j$. Let $f_I=\sum_{i \in I}p_if$ and $f_J=\sum_{j \in J}p_j f$. It is not difficult to check that $E= E_I\oplus E_J$, $[{E_I}_\lambda E_J]=0$, $f= f_I+ f_J$, where $f_I$ is a homomorphism from $\mathcal{A}$ to $E_I$ and $f_J$ is an anti-homomorphism from $\mathcal{A}$ to $E_J$. This completes the proof.
\end{proof}

{\bf Acknowledgments}
{This work was supported by the National Natural Science Foundation of China (Nos. 12171129, 12271406) and the Fundamental Research Funds for the  Central Universities (No. 22120210554).}

\end{document}